\documentclass[11pt,a4paper]{article}
\usepackage[a4paper]{geometry}
\usepackage{amssymb,latexsym,amsmath,amsfonts,amsthm}
\usepackage{graphicx}
\usepackage{epstopdf}
\usepackage{tikz}
\usepackage{pgflibraryshapes}
\usetikzlibrary{arrows,decorations.markings}
\usepackage{comment,verbatim}
\usepackage{hyperref}
\newtheorem{thm}{Theorem}[section]
\newtheorem{lem}{Lemma}[section]
\newtheorem{pro}{Proposition}[section]

\newtheorem{exa}{Example}[section]

 \newtheorem{rem}{Remark}[section]

 \theoremstyle{remark}

\numberwithin{equation}{section}

\def\Im {\mathop{\rm Im}\nolimits}
\def\arg {\mathop{\rm arg}\nolimits}
\def\Re {\mathop{\rm Re}\nolimits}

\def\Ai {{\rm Ai}}
\def\Bi {{\rm Bi}}
\def\arccosh {{\rm arccosh}}

\numberwithin{equation}{section}

\begin{document}

\title{ Recent Advances in Asymptotic Analysis}

\author{
 R. Wong\footnotemark [0]\footnotemark [2]  ~and Yu-Qiu Zhao\footnotemark [3] }

\renewcommand{\thefootnote}{\fnsymbol{footnote}}

\footnotetext [0]   {{\bf{Funding:}} The second author was  supported in part by the National Natural  Science Foundation of China under grant numbers    11571375 and 11971489. }

\footnotetext [2]  {Department of Mathematics, City University of Hong Kong, Tat Chee Avenue,
Kowloon, Hong Kong.}

\footnotetext [3]  { Department of Mathematics, Sun Yat-sen University, GuangZhou 510275, China.}

\date{}
\maketitle

 \begin{abstract}
 This is a survey article on an old topic in classical analysis. We present some new developments in asymptotics in the last fifty years. We start with the classical method of Darboux and its generalizations, including an uniformity treatment which has a direct application to the Heisenberg polynomials. We then present the development of an asymptotic theory for difference equations, which is a major advancement since  the work of
  Birkhoff and   Trjitzinsky
  in 1933. A new method was introduced into this field in the nineteen nineties, which is now known as the nonlinear steepest descent method or the Riemann-Hilbert approach. The advantage of this method is that it can be applied to orthogonal polynomials which do not satisfy any differential or difference equations neither do they have any integral representations. As an example, we mention the case of orthogonal polynomials with respect to the Freud weight. Finally, we show how the Wiener-Hopf technique can be used to derive asymptotic expansions for the solutions of an integral equation on a half line.


 \end{abstract}

\tableofcontents

\section{Introduction}\label{sec: introduction}
Asymptotics is one of the oldest topics in classical analysis; see \cite{Birkhoff-1973}.
Examples from this topic include Stirling's formula \cite{Dutka-1991,Stirling-1730}, prime number theorem \cite{Titchmarsh-1951}, and
Hardy-Ramanujan's formula \cite{Hardy-Ramanujan-1918}. 
 Despite of the fact that it is an old topic, new methods continue to be discovered. For instance, in the construction of error bounds for asymptotic expansions of the Stieltjes transform and the fractional integral transform, a method was introduced which makes use of distribution theory; see \cite{McClure-Wong-1978} and \cite{McClure-Wong-1979}. This method was later extended to derive asymptotic expansions of the generalized  Mellin transform \cite{Wong-McClure-1984}, which is a divergent integral but exists in the sense of generalized functions. In a study of the behavior of polynomials orthogonal with respect to the Freud weight $\exp\left(-|x|^\beta\right)$, $\beta>0$, another new method was brought into the field of asymptotics, which is the Riemann-Hilbert approach  introduced by Deift-Zhou \cite{Deift-Zhou-1993} in their investigation of the behavior of solutions to the mKdV equation.

The development of an asymptotic theory for ordinary differential equations  was quite successful in the last century. After Poincar\'{e} \cite{Poincare-1886} introduced a proper definition for asymptotic expansion and established rigorously the results for irregular singularity at infinity, several people obtained asymptotic expansions for solutions of differential equations with a large parameter, including   Schlesinger \cite{Schlesinger-1907}, Birkhoff \cite{Birkhoff-1908}, and  Turrittin \cite{Turrittin-1936}. More important results were given by Langer \cite{Langer-1931,Langer-1932,Langer-1949} on uniform asymptotic solutions to turning-point problems, and his results were substantially extended by Olver
\cite{Olver-1954a,Olver-1954b, Olver-1956} including extensions to unbounded domains of the independent  variable \cite{Olver-1956}
and explicit bounds for the error terms \cite{Olver-1961}.
The corresponding development for difference equations, however, was not as successful as for differential equations. Formal theory of irregular linear difference equations was first given by Birkhoff \cite{Birkhoff-1930}, and rigorous analysis was later given by Birkhoff and Trjitzinsky \cite{Birkhoff-Trjitzinsky-1933}.  The first paper appeared in 1930, and the second paper was published in 1933.  For the remaining part of the century, there were essentially no significant results in this area of research.  The first paper on turning-point problem for difference equation \cite{Wang-Wong-2002} did not appear until 2002; see also \cite{Wong-2014}. Subsequently, progress seems to have picked up, and results in this area are now nearly as complete as they are for differential equations; see \cite{Wang-Wong-2003},  \cite{Wang-Wong-2005}, \cite{Cao-Li-2014} and \cite{Li-Wang-Wong-2020}.

While there are now plenty of results on finding asymptotic solutions to differential equations and to difference equations, the same is not true for integral equations. The problem is not that difficult if the equation is either of Volterra type or Fredholm type with a difference kernel, since the equation can be reduced to an algebraic equation by taking either Laplace transforms or Fourier  transforms, and asymptotic behaviour of the solution can be obtained by taking their inverse transformations; see, e.g., \cite{Wong-Wong-1974}, \cite[p.\;286]{Beals-Wong-2020}. The situation is quite different if the equation is of the form
\begin{equation}\label{int-eq}
u(t)=f(t)+\int^\infty_0 k(t-\tau) u(\tau) d\tau.
\end{equation}
As far as we are aware, there is only one paper on asymptotic solutions to such an equation \cite{Muki-Sternberg-1970}. Equation in \eqref{int-eq} is known as of Wiener-Hopf type or renewal equation on half line. The famous Wiener-Hopf technique is to construct a contour-integral representation for the solution of \eqref{int-eq}, but to find the behaviour of the solution is another story since the functions involved in the integrand are very complicated and contain Stieltjes  and Hilbert transforms. A rather satisfactory solution to this problem was given only very recently; see \cite{Li-Wong-2021}.

The purpose of this survey is to bring people's attention to some of the developments in asymptotic analysis in the past 50 years. In Section\;\ref{sec: Darboux}, we begin with Darboux's method. Even this classical result has some very attractive extensions. We will first present a generalization that involves logarithmic singularities \cite{Wong-Wyman-1974}, and mention couple examples to show that this generalization is indeed useful. Next, we give a result on the behavior of the coefficients in a Maclaurin series, whose generating function has two coalescing algebraic singularities on its circle of convergence, and provide as an application the asymptotic expansion of the Heisenberg polynomials as the degree of the polynomial goes to infinity; see \cite{Liu-Wong-Zhao-2013} and \cite{Wong-Zhao-2005}.

In Section\;\ref{sec: difference}, we briefly review the results that have been obtained for difference equations during the past twenty years. As examples, we present the asymptotic expansions of the Racah and the Wilson polynomials; these polynomials are listed at the top level of the Askey scheme \cite[p.\;464]{nist-2010}.

In Section\;\ref{sec: Riemann-Hilbert}, we introduce the Riemann-Hilbert approach, and show how to use it to derive the asymptotic behavior of the polynomials orthogonal with respect to the Freud weight $\exp\left (-|x|^\beta\right )$, $\beta>0$. If $\beta\geq 1$, then the problem was completely solved by
Kriecherbauer and McLaughlin \cite{Kriecherbauer-McLaughlin-1999} in 1999. The difficulty of this problem occurs when  $0<\beta<1$, in which case it is not known what is the dominant approximant in the asymptotic expansion of these polynomials; see \cite[p.\;62]{Deift-Kriecherbauer-McLaughlin-Venakides-Zhou-2001}. This problem was finally settled in 2016, and the answer is that the dominant approximant does not involve special functions such as Airy or Bessel; it involves two  solutions to a pair of decoupled scalar  linear integral equations.

In the final section, Section\;\ref{sec: integral-eq}, we return to the integral equation  \eqref{int-eq}. We first recall what is the Wiener-Hopf technique \cite{Wiener-Hopf-1931}, the solution to equation \eqref{int-eq} is then given by a contour integral.  Thus,  we  reduce the problem to an asymptotic evaluation of a Fourier integral whose integrand involves Stieltjes and Hilbert transforms.


\section{Darboux's method}\label{sec: Darboux}

The following  combinatorial problem has appeared  in Whitworth \cite{Whitworth-1901}.

{\it If there
be $n$ straight lines in one plane, no three of which meet in a point, no two lines are parallel, there will be $n(n-1)/2$ points of intersection.
Put the intersection points into groups of $n$, in each of which no three
points
lie in one of the straight lines.
Find  the number $g_n$ of such groups.}

 The answer given in \cite{Whitworth-1901} is  $g_n= \frac 1 2(n -1)!$ for $n=3,4,\cdots$.  However,
in 1951, Robinson \cite{Robinson-1951} considered the problem and showed that the answer is wrong.
He found the correspondence between $n$-point groups and $n$-polygons.
By counting  $g_n=\frac 1 2(n -1)!$, one misses the composite $n$-gons:
those  consist  of several separate simple polygons with totally $n$ sides.
Instead, $u_n=g_n/\left(\frac 1 2 (n-1)!\right)$ satisfies the recurrence
relation
\begin{equation}\label{Robinson-DE}
u_{n+3}=u_{n+2}+\frac 1 { 2 n } u_{n},~~~n= 1,2,3,\cdots,
\end{equation}with  $u_1=u_2=0$ and $u_3=1$. After confirming the existence of   the number
\begin{equation*}
b=\lim_{n\to\infty}\frac{u_n^2} n=0.284098\cdots,
\end{equation*}  as noted by   Robinson \cite{Robinson-1951},
 $b$  assumes the role of an absolute geometric constant. Robinson further asked  whether it is an algebraic number, which seems unlikely,
whether it is expressible in terms of $\pi$ and $e$, or whether it is entirely new.

The answers came fast.
It is stated that the constant is given by
\begin{equation*}
b=4 e^{-3/2}/\pi;
\end{equation*} cf. an
editorial note \cite{Robinson-constant-1952}.
Thirteen readers, including 
 Harry Pollard,  John Riordan,
George Szekeres, Fritz Ursell  and Morgan Ward, wrote to the editor about Robinson's constant.   Several of them     used the generating function
\begin{equation}\label{Robinson-Ge}
f(x)=-2+2(1-x)^{-  1/ 2} \exp\left( - \frac{x^2+2x} 4\right )=\sum^\infty_{n=3} \frac{u_n} n x^{n},
\end{equation}which can readily be derived from \eqref{Robinson-DE} via the differential equation
\begin{equation*}
2(1-x)f'(x)=x^2\left [f(x) +   2\right] ,~~~~f(0)=0.
\end{equation*} Mention had also been made of a theorem from  Titchmarsh \cite[p.\;224]{Titchmarsh-1958}:
\begin{pro}\label{prop:Titchmarsh}
Let
 \begin{equation*}
 f(x)=\sum^\infty_{n=0} a_n x^n,~~~~g(x)=\sum^\infty_{n=0} b_n x^n,~~~~|x|<1,
 \end{equation*}
  where $a_n\geq 0$, $b_n\geq 0$,   and
$\displaystyle{\sum^\infty_{n=0} b_n }$ diverges.  If  $\displaystyle{\lim_{n\to \infty } a_n/b_n=C }$, then
$\displaystyle{\lim_{x\to 1 } {f(x)}/{g(x)} =C }$.
\end{pro}

To determine the behavior of $a_n$, it seems that  the converse of   Proposition \ref{prop:Titchmarsh} matters much.   Yet   Titchmarsh \cite[p.\;226]{Titchmarsh-1958}
claimed that there is no general converse of the above proposition:   {\it from the asymptotic behavior of $f(x)$ we can not deduce that of $a_n$}.

Nevertheless,  Darboux's method (G. Darboux, 1878, \cite{Darboux-1878}) demonstrates how to  deduce the asymptotic behavior of $a_n$   from that of $f(x)$.
Such problems arise in various mathematical fields, including
 number theory,   combinatorics, and orthogonal polynomials.
 As a start, one may use the Cauchy integral formula
 \begin{equation*}
  a_n=\frac 1 {2\pi i} \int_C
f(z) z^{-n-1} dz ,
 \end{equation*}
where $C$ is a circle encircling the origin, but not any singularities of $f(z)$.
Expanding the contour $C$ if possible,
and applying the residue calculus if  a   pole  is encountered, we obtain the leading contribution to the asymptotic behavior of $a_n$.
This procedure fails, however, if an
  algebraic singularity is the closest to the origin, and  the classic Darboux's method (\cite{Darboux-1878}; see also \cite{Szego-1975} and \cite{Wong-1989}) becomes   useful.

To describe briefly how
  this method works,
we assume that
\begin{equation}\label{generating-func}
f(z)=\sum^\infty_{n=0} a_n z^n
\end{equation}
is analytic in  $|z|<1$, with only one singularity on the circle $|z|=1$, say, at $z=1$, and that in a neighborhood of $z=1$, $f(z)$ is of the form
\begin{equation}\label{f-g-Darboux}
f(z)=(1-z)^\alpha g(z),
\end{equation}
 where $g(z)$ is analytic at $z=1$, and $\alpha\not\in \mathbb{Z}$.
 To find the large-$n$ behavior of $a_n$, we expand
  \begin{equation}\label{g-Darboux}
  g(z)=\sum^\infty_{r=0} c_r(1-z)^r
  \end{equation}at $z=1$.
   The $m$-th Darboux approximant of $f(x)$ is then defined by
  \begin{equation*}
 f_m(z) =\sum^m_{r=0} c_r(1-z)^{\alpha +r},
  \end{equation*}which is also analytic in $|z|<1$; cf. \eqref{generating-func}.
  Expanding $f_m(z)$ into a Maclaurin series
 \begin{equation*}
  f_m(z)=\sum^\infty_{n=0} b_{mn}z^n,
 \end{equation*}
   it turns  out $a_n\sim b_{mn}$ for finite $m$, as $n\to \infty$.
Indeed, we can write
 $a_n= b_{m,n}+\varepsilon_m(n)$, and show that  the error term
\begin{equation*}
 \varepsilon_m(n)=\frac 1 {2\pi i}
\int_C \left [f(z)-f_m(z) \right] z^{-n-1} dz =o(n^{-N}),
\end{equation*} where  $\Re\alpha+m+1\leq N< \Re\alpha +m+2.$
The above error estimate has been proved rigorously; see, e.g., Wong \cite[Ch.\;2.6]{Wong-1989}. Thus we have the following   result.
\begin{thm}\label{Thm:Darboux-classic-1-p}
With $f(z)$ given in \eqref{generating-func}, \eqref{f-g-Darboux} and \eqref{g-Darboux}, assume that $\alpha\not\in \mathbb{Z}$ is a complex constant,
and $z=1$ is the only singularity on the circle of convergence $|z|=1$. For any $m\geq 0$, the Maclaurin coefficient of $f(z)$ has the  asymptotic expansion
\begin{equation}\label{expansion-Darboux}
a_n=\sum^m_{r=0} c_r  \begin{pmatrix}
                                       \alpha+r\\
                                       n \\
                                       \end{pmatrix}(-1)^n
 +o\left(n^{-\Re\alpha-m-1}\right )
\end{equation}as
$n\to\infty.$
\end{thm}

Formula \eqref{expansion-Darboux} is referred to as a
generalized asymptotic expansion in \cite{Wong-1989}, which is more general than the usual Poincar\'{e} expansion; see Erd\'{e}ly and    Wyman \cite{Erdely-Wyman-1963}.
The asymptotic nature of the expansion is rooted in
\begin{equation*}
 \begin{pmatrix}
\alpha+r\\
n \\
\end{pmatrix}(-1)^n=\frac {\Gamma(n-\alpha-r)}{n!\Gamma(-\alpha-r)}
\sim \frac {  n^{-\alpha-1-r}}{ \Gamma(-\alpha-r)}
 , ~~~r=0,1, \cdots~~\mbox{as}~n\to\infty .
\end{equation*}

When there are finite number of algebraic singularities on the circle of convergence, the final result can be extended by
adding  up the contributions from all singular points; see Szeg\H{o} \cite[p.\;207, Thm.\;8.4]{Szego-1975}  and Wong \cite[p.\;119, Thm.\;5]{Wong-1989}.

Applying Theorem \ref{Thm:Darboux-classic-1-p} to the generating function \eqref{Robinson-Ge} (cf. \cite[p.\;120]{Wong-1989}), one finds that $u_n\sim b\sqrt n$ as $n\to\infty$, with $b=4 e^{-3/2}/\pi$; this is Robinson's constant.
The theorem can also be used to solve a problem
of
 Erlebach and  Ruehr
 in counting Hamiltonian
cycles for bipartite graphs; see
 Knuth \cite{Knuth-1980}.
Darboux's method plays an important role as well in determining    orthogonal measures, especially those with explicit three-term recurrence relations; see
Ismail \cite{Ismail-2005}.

We complete our discussion of the classical Darboux's method by mentioning the following interesting example investigated by Olver \cite{Olver-1970}.
It is well known that the Legendre polynomials possess   the generating function
\begin{equation}\label{Legendre-g}
f(z)=\left [(e^{i\theta}-z)(e^{-i\theta}-z)\right ]^{-1/2}=
\sum^\infty_{n=0} P_n(\cos\theta) z^n,
\end{equation}
where the branches are chosen such that  $(e^{\pm i\theta}-z)^{-\frac
1 2}\rightarrow e^{\mp \frac 1 2 i\theta}$ as $z\rightarrow 0$.
Now applying
  Darboux's method yields
\begin{equation*}
P_n(\cos\theta) \sim \left (\frac 2 {\sin\theta}\right )^{\frac
1 2} \sum^\infty_{k=0}\left (
\begin{array}{c}
-\frac 1 2 \\k\end{array} \right )\left (
\begin{array}{c}
k-\frac 1 2 \\n\end{array} \right )\frac
{\cos\theta_{n,k}}{(2\sin\theta)^k},~~~~n\to\infty,
\end{equation*}
 where
\begin{equation*}
 \theta_{n,k}
=\left(n-k+\frac 1 2\right)\theta +\left(n-\frac 1 2 k-\frac 1 4\right)\pi.
\end{equation*}
On the other hand,
\begin{equation*}
P_n(\cos\theta) = \left (\frac1  {2\sin\theta}\right )^{\frac
1 2} \sum^\infty_{k=0}\left (
\begin{array}{c}
-\frac 1 2 \\k\end{array} \right )\left (
\begin{array}{c}
k-\frac 1 2 \\n\end{array} \right )\frac
{\cos\theta_{n,k}}{(2\sin\theta)^k}
\end{equation*}
    for $\frac 1 6
\pi<\theta<\frac 56 \pi$.
A paradox now arises, since one has
 \begin{equation}\label{paradox}
 P_n(\cos\theta) \sim 2P_n(\cos\theta ), ~~~~~\frac 1 6
\pi<\theta<\frac 56 \pi,~~~~~ n\to\infty ;
 \end{equation}cf. Olver \cite{Olver-1970}.

 \vskip 1cm

Darboux's method dates back to 1878. Although frequently used, it does not appear to
have been extended until around 1970.
In 1974, Wong and Wyman \cite{Wong-Wyman-1974}  have
given a generalization of Darboux's method, which allows the generating function $f(z)$
in \eqref{generating-func} to have logarithmic-type singularities on its circle of convergence. This paper was considered {\it an important early (and unduly neglected) reference} in analytic combinatorics;
see Flajolet and Sedgewick \cite[p.\;438]{Flajolet-Sedgewick-2009}.
First we quote a couple of examples of this  type.

{\bf{Problem of P\'{o}lya}} (1954, \cite[Ex.\;8,\;p.\;8 and p.\;213]{Polya-1954}): Define $A_n$ by
 \begin{equation}\label{Polya}
 \frac z {\ln(1+z)}   =\sum^\infty_{n=0} A_n \frac{z^n}{n!},~~~~~~|z|<1.
 \end{equation}Show that
 \begin{equation}\label{Polya-approx}
 \frac{A_n }{n!}\sim (-1)^{n-1}\frac {1} {n \ln^2 n}, ~~\mbox{as}~n\to \infty.
 \end{equation}

{\bf{Problem of Knuth}} (2004, \cite{Knuth-2004}): What is the asymptotic behavior of the constant $l_n$ defined by the generating function
\begin{equation}\label{Knuth-2004}
 \frac 1 {\ln (1-z)}+\frac 1 z=\sum^\infty_{n=0} l_n z^n,~~~~~~|z|<1.
 \end{equation}

 In \cite{Wong-Wyman-1974},  Wong and  Wyman dealt with the general case
\begin{equation}\label{Wong-Wyman}
f(z)=(1-z)^\alpha [ \ln (1-z) ]^\mu g(z)=\sum^\infty_{n=0} a_n z^n, ~~~~|z|<1,
\end{equation}
 where $\alpha$ and $\mu$ are complex constants, and
$g(z)$ is analytic in $|z|\leq 1+\delta$ for some $\delta>0$. Now we briefly outline their treatment. Again the Cauchy integral formula applies, and one has, by deforming the contour,
\begin{equation*}
 2\pi i a_n= \oint_{|z|=1+\delta_n} f(z) z^{-n-1} dz
-\oint_{|z-1|=\delta_n} f(z) z^{-n-1}
dz,
\end{equation*}
 where the integration paths are counterclockwisely  oriented, and  $\delta_n=n^{-  1 / 2}$. The first integral is of the order   $O ( e^{-c\sqrt n} \; )$ for some positive constant $c$; the second integral alone contributes to the full asymptotic expansion.
The approximants are   ``special'' functions with  integral representations
\begin{equation}\label{special-function-WW}
 M(\alpha, \mu, n)=-\frac 1 {2\pi i} \int^{(0+)}_{+\infty} (-t)^\alpha
[\ln(-t)]^\mu e^{-(n+1)t} dt,~~~n=0,1,2,\cdots,
\end{equation}and have full asymptotic expansions of the form
\begin{equation*}
M(\alpha, \mu, n)\sim \frac {(-\ln(n+1))^\mu }{(n+1)^{\alpha+1}}
\sum_{k=0}^\infty \begin{pmatrix}
                    \mu \\
                    k \\
                  \end{pmatrix}
 \frac {d^k}{d\alpha^k}\left[ \frac 1 {\Gamma(-\alpha)} \right ] \frac 1  {(-\ln(n+1))^k }.
\end{equation*}
Wong and Wyman further introduced the function
\begin{equation}\label{J-WW}
J_m=-\frac 1 {2\pi i} \oint_{|t|=\delta_n}
(-t)^{\alpha+m}
[\ln(-t)]^\mu P_m((n+1)t) e^{-(n+1)t} dt,
\end{equation}where $P_m(w)$ is a polynomial given by
\begin{equation*}
P_m(w)=\left . \frac 1{m!} \frac {d^m}{dt^m} \left [
g(t+1)\exp\left\{ -\frac 1 2 wt \left [\frac{2(\ln(1+t)-t)}{t^2}\right ]\right\}\right ]\right |_{t=0}:=\sum^m_{s=0} p_sw^s,
\end{equation*}so that
\begin{equation*}
J_m\sim \sum^m_{s=0} p_s (n+1)^s M(\alpha+m+s, \mu, n)
 ~~~~\mbox{as}~n\to\infty.
\end{equation*}

\begin{thm}\label{thm:Wong-Wyman}(Wong-Wyman, \cite{Wong-Wyman-1974})
If $f(z)$ is the function given in \eqref{Wong-Wyman}, then its Maclaurin coefficient $a_n$ has the asymptotic expansion
\begin{equation*}
a_n\sim \sum^\infty_{m=0} (-1)^m J_m,
\end{equation*}where $J_m$ is defined in \eqref{J-WW} and has the asymptotic expansion
\begin{equation*}
J_m\sim     \frac {\left(-\ln (n+1)\right )^\mu}{(n+1)^{\alpha+m+1}}
\sum_{k=0}^\infty \left\{\begin{pmatrix}
                    \mu \\
                    k \\
                  \end{pmatrix}
\sum^m_{s=0} p_s   \frac {d^k}{d\alpha^k}\left[ \frac 1 {\Gamma(-\alpha-m-s)} \right ]\right\} \frac 1  {(-\ln(n+1))^k }
,~n\to\infty.
\end{equation*}
\end{thm}

\vskip .5cm

The above result includes Theorem \ref{Thm:Darboux-classic-1-p} as a special case if $\mu=0$, and  can be extended to cases with  several fixed  singularities on the circle of convergence.

Now it is readily seen that
\eqref{Polya-approx} follows accordingly, by applying Theorem \ref{thm:Wong-Wyman} to \eqref{Polya}.
For the generating function  \eqref{Knuth-2004}, we take
 $\alpha=0$, $\mu=-1$ and $g(z)\equiv 1$ in \eqref{Wong-Wyman}. Then, Theorem \ref{thm:Wong-Wyman}
applies  and we have
 \begin{equation*}
 l_n\sim \frac 1 {(n+1) \ln^2(n+1)} ~~~~\mbox{as}~n\to\infty.
 \end{equation*}As a consistence check, it is worth noting that \eqref{Polya-approx} and \eqref{Knuth-2004} are related in the straightforward manner: $A_n/n!=(-1)^{n-1} l_{n-1}$ for positive integers $n$.
\vskip .5cm

When the singularities are free to move on the circle of convergence, Darboux's method
will continue to work only if their essential configuration remains the same while the relative
positions vary. However, this method breaks down when two or more singularities coalesce with each other.
For example, the generating function
\eqref{Legendre-g}
for the Legendre polynomials
has a pair of algebraic singularities $e^{\pm i\theta}$. Darboux's method fails to apply when  $e^{\pm i\theta}\to 1$ as
$\theta\to 0^+$, and $e^{\pm i\theta}\to -1$ as
$\theta\to \pi^-$; a difficulty arises.
In the  cases of   coalescing singular points, the asymptotic expansion may involve transcendental functions
instead of elementary ones.

In 1967, Fields \cite{Fields-1967} presented  a uniform treatment of Darboux's method when two or three
singularities coalesce. More precisely, he considered the case in which
\begin{equation}\label{Fields}
f(z,\theta) = (1-z)^{-\lambda} \left [ \left(e^{ i\theta}-z\right)  \left(e^{-i\theta}-z\right)\right ]^{-\Delta}
g (z,\theta)=
\sum^\infty_{n=0}a_n(\theta) z^n,
\end{equation}
where the Maclaurin expansion converges for $|z|<1$, $\lambda$ and $\Delta$ are bounded quantities,
the branches of
$(1-z)^{-\lambda}$ and  $\left [ \left(e^{ i\theta}-z\right)  \left(e^{-i\theta}-z\right)\right ]^{-\Delta}$
  are chosen such that each of
them reduces to $1$ at $z = 0$, and $g(z, \theta)$ is analytic in $|z|\leq e^\eta$ ($\eta >0$), uniformly for
$\theta\in [0,\pi]$.
  In \cite{Fields-1967}, Fields first expressed $a_n(\theta)$   as a Cauchy integral, then made a change of
variable and a rescaling, and finally  obtained a generalized asymptotic expansion in the
sense of Erd\'{e}lyi and Wyman \cite{Erdely-Wyman-1963}. His results are uniform in certain $\theta$-intervals depending
on $n$.

Despite the fact that Fields' results have achieved the so-called uniform reduction in
the sense of Olver \cite[p.\;102]{Olver-1975}, they are found to be too complicated for any practical
application; see, e.g., Erd\'{e}lyi \cite[p.\;167]{Erdely-1970}, Olver \cite[pp.\;112-113]{Olver-1975}, and Wong \cite[p.\;145]{Wong-1989}.
For example, Olver commented that
 {\it{It may be desirable to investigate
whether any simplifications are feasible since the results in \cite{Fields-1967} are rather complicated
to apply in their present form}.}

 A   motivation for the uniform treatment of Darboux's method  came from  the uniform asymptotic expansions of integrals. As a point of information, we mention a few  relevant references in this respect. For example,
 Chester, Friedman, and
Ursell (1957, \cite{Chester-Friedman-Ursell-1957}) first  presented a uniform treatment  of the steepest descent method,
Bleistein (1967, \cite{Bleistein-1967}) considered the problem of many nearby stationary points and algebraic singularities,
Wong (1989, \cite[Ch.\;VII]{Wong-1989}) provided various types of coalescence of singular points, in particular uniform asymptotics of orthogonal polynomials, and
Olde Daalhuis and Temme (1994, \cite{Olde-Daalhuis-Temme-1994}) introduced  a class of  rational functions based on which
  error terms can be estimated.

The work of  Wong and  Zhao  \cite{Wong-Zhao-2005} in 2005
seems to have responded to Olver's request mentioned above,
   i.e., to derive simpler forms of uniform asymptotic expansions when two
or more algebraic singularities, on the circle of convergence, coalesce with each other
as  certain parameter approaches a critical value.
To begin with,
Wong and  Zhao first considered the   simplest case
with two singularities, namely
\begin{equation}\label{Wong-Zhao}
 f(z,\theta) =\left [ (e^{i\theta}-z)(e^{-i\theta}-z)\right ]^{-\alpha}g(z, \theta ) = \sum^\infty_{n=0} a_n(\theta) z^n,
\end{equation}where
$g(z, \theta)$ is analytic in $|z|\leq e^\eta$, $\eta>0$, which is a special case of Fields; cf. \eqref{Fields} with $\lambda=0$.
Contribution to the large-$n$ behavior still comes from the singular points $e^{\pm i\theta}$, which are now  subject to vary as $\theta\to 0^+$. The approximants in this case are
\begin{equation*}
T_1(x) = \frac 1{2\pi i}\int_{\Gamma_0}(s^2+1)^{-\alpha}   e^{xs} ds,~~~~ T_2(x) = \frac 1{2\pi i}\int_{\Gamma_0}s(s^2+1)^{-\alpha}   e^{xs} ds,
\end{equation*}
where $\Gamma_0$ is a Hankel-type loop which starts and ends at $-\infty$, and encircles $s=\pm i$ in
the positive sense.
An asymptotic expansion of the form
\begin{equation*}
 a_n(\theta)=  \theta^{1-2\alpha}T_1(n\theta)\sum^{m-1}_{k=0}\frac {\alpha_k(\theta)}{n^k} + \theta^{1-2\alpha} T_2(n\theta)\sum^{m-1}_{k=0}\frac {\beta_k(\theta)}{n^k}  + \varepsilon(\theta, m)
\end{equation*}was obtained for
 $m=1,2,3,\cdots$, with the error term
 $\varepsilon(\theta, m)$
 given explicitly and estimated, and the coefficients determined recursively as follows:
Begin with an analytic function
\begin{equation}\label{h-0}
h_0(s,\theta) :=g(e^{-s\theta}, \theta)\left [\left (\frac{e^{-s\theta}-e^{i\theta}}{(-s-i)\theta}\right )\left(\frac{e^{-s\theta}-e^{-i\theta}}{(-s+i)\theta}\right )\right ]^{-\alpha}
\end{equation}
  for $\Re s\geq -\eta/\theta$ and $|s\pm i|<2\pi/\theta$,
and make use of the
 iteration
\begin{equation}\label{iteration}
\left\{ \begin{array}{l}
  \displaystyle{ h_k(s,\theta):=\alpha_k(\theta ) +s {\beta_k(\theta )} +
\left (s^2 +1 \right )g_k(s,\theta ),}    \\[.4cm]
\displaystyle{   h_{k+1}(s,\theta)
=-
\frac 1\theta\left [ (s^2+1) \frac d {ds} +2(1-\alpha )s \right ]g_k(s,\theta),  }
        \end{array}
\right .
k=0,1,2,\cdots .
\end{equation}
Assuming analyticity of each $h_k(s,\theta)$ and $g_k(s,\theta)$ at $s=\pm i$,
one determines the coefficients $\alpha_k(\theta )$ and $\beta_k(\theta )$.
 Furthermore, straightforward verification shows that
\begin{equation*}
T_1(x)=\frac{\sqrt\pi} {\Gamma(\alpha)} \left (\frac {x} 2\right )^{\alpha-\frac 1 2}
J_{\alpha-\frac 1 2}(x),
\end{equation*}and
\begin{equation*}
T_2(x)=\frac{\sqrt\pi} {\Gamma(\alpha)} \left (\frac {x} 2\right
)^{\alpha-\frac 1 2}\left [ \frac{2\alpha-1  } {x}  J_{\alpha-\frac 1 2}(x)
- J_{\alpha+\frac 1 2}(x)\right ]=\frac{\sqrt\pi} {\Gamma(\alpha)} \left (\frac {x} 2\right
)^{\alpha-\frac 1 2} J_{\alpha-\frac 3 2}(x).
\end{equation*}
Hence we have

\begin{thm}(Wong-Zhao \cite{Wong-Zhao-2005}; see also  \cite[Cor.\;1]{Liu-Wong-Zhao-2013})
If $f(z,\theta)$ is the function given in \eqref{Wong-Zhao}. Then its Maclaurin coefficient has the asymptotic expansion
\begin{equation}\label{Wong-Zhao-expansion}
a_n(\theta)
\sim \frac{\sqrt\pi} {\Gamma(\alpha)}\left (\frac {n}{2\theta}\right )^{\alpha-\frac 1 2}\left [ J_{\alpha-\frac 1 2}(n\theta) \sum^{\infty}_{k=0}\frac { \alpha_k(\theta)}{n^k} + J_{\alpha-\frac 3 2}(n\theta)  \sum^{\infty}_{k=0}\frac { \beta_k(\theta)}{n^k}\right ]
\end{equation}as $n\to \infty$, holding uniformly for $\theta \in [0, \pi-\delta]$, $\delta>0$, with coefficients
  $\alpha_k(\theta )$ and $\beta_k(\theta )$  given by \eqref{h-0} and \eqref{iteration}. The leading coefficients are
  $\alpha_0(\theta)=\cos(\alpha\theta)\left(\frac{\sin\theta}{\theta}\right )^{-\alpha}$ and $\beta_0(\theta)=\sin(\alpha\theta)\left(\frac{\sin\theta}{\theta}\right )^{-\alpha}$.
\end{thm}\vskip .5cm

The above result can be extended to  handle cases with many coalescing algebraic singularities.  For instance, we have
\begin{equation}\label{Wong-Zhao-many}
 f(z,\theta)=\left\{ \prod^q_{k=1} (z_k(\theta)-z)^{-\alpha_k}\right\} g(z,\theta) =
\sum^\infty_{n=0} a_n(\theta ) z^n,
\end{equation}where
$g(z, \theta)$ is analytic in $|z|\leq e^\eta$ with $\eta>0$,
$z_k(\theta)\equiv e^{i \theta s_k(\theta)}$ so that $z_k(\theta)\rightarrow 1$  as $\theta\rightarrow 0$.
In this case,  the approximants are the special functions
\begin{equation}\label{T-L-def}
 T_l(x):=
\frac 1 {2\pi i}\int_{\Gamma_0}
s^{l-1}e^{x s}
 \prod^q_{k=1} (s+is_k(\theta))^{-\alpha_k}  ds
\end{equation}for
$l=1,2,\cdots, q$, where $\Gamma_0$ is also a Hankel-type loop which starts and ends at $-\infty$, encircling all points $s=-is_k(\theta)$ in the positive direction.

\begin{thm}(Wong-Zhao \cite{Wong-Zhao-2005})\label{thm:Wong-Zhao-many}
Let $f(z,\theta)$ be the function given in \eqref{Wong-Zhao-many}. If
$z_k(\theta)=e^{is_k\theta}$, $s_k$ being real constants, then
\begin{equation}\label{Wong-Zhao-many-expansion}
 a_n(\theta)\sim \theta^{1-\alpha}\sum^{q}_{l=1}T_l(n\theta)\sum^{\infty}_{k=0}\frac {\beta_{k,l}(\theta)}{n^k}
\end{equation}as $n\to \infty$, uniformly in $\theta\in [0, \nu]$ with $\nu<
\min_{1\leq k\leq q} \{ \pi/ |s_k|\}$, where
  ${\alpha=\sum^q_{k=1}\alpha_k}$, and the coefficients $\beta_{k,l}(\theta)$ can be determined successively.
\end{thm}\vskip .5cm

Properties of the special functions  $T_l(x)$ in \eqref{T-L-def} with $s_k(\theta)\equiv s_k$  are considered in \cite{Wong-Zhao-2005}.
For instance, from \eqref{T-L-def} we have
\begin{equation*}
 T_l(x)=\frac {d^{l-1}}{d x^{l-1}} T_1(x)~~~~\mbox{for}~~l=1,2,\cdots, q ,
\end{equation*}and that $w=T_1(x)$ solves the differential equation
 \begin{equation*}
 x\frac {d^{q}w}{d x^{q}}  +\sum^{q-1}_{l=0}
(d_l x +c_l)\frac {d^{l }w}{d x^{l }}  =0,
 \end{equation*}
where the coefficients $c_l$ and $d_l$ are determined by
\begin{equation*}
\prod_{k=1}^q (s+is_k):=s^q+\sum^{q-1}_{l=0} d_l s^l,~~~~
 \sum^q_{k=1} (1-\alpha_k)
\prod_{k'\not=k} (s+is_{k'}):=\sum^{q-1}_{l=0} c_l s^l.
\end{equation*}

Now we check the special case \eqref{Fields} investigated by
 Fields \cite{Fields-1967}. Theorem \ref{thm:Wong-Zhao-many}
 applies,  and we may set
\begin{equation*}
 q=3;~~~~s_1=0,~\alpha_1=\lambda;~~~~s_2=1,~\alpha_2=\Delta;~~~~
s_3=-1,~\alpha_3=\Delta.
\end{equation*}By expanding a part of the integrand in \eqref{T-L-def}
into   Laurent series,  we find that
 in terms of the generalized hypergeometric functions, we have
\begin{equation*}
 T_l(x)= \frac{x^{\lambda+2\Delta-l}}{\Gamma(\lambda+2\Delta-l+1)}\,
 {}_1 F_2\left (\Delta; \frac {\lambda+2\Delta-l+1} 2, \frac {\lambda+2\Delta-l+2} 2;-\frac {x^2} 4\right )
\end{equation*}
for $l=1,2,$ and $3$; cf. \cite[(16.2.1)]{nist-2010}.

We complete our discussion on uniform treatment of Darboux's method by listing several follow-up progresses.
 In 2007, Bai and  Zhao \cite{Bai-Zhao-2007} derived uniform asymptotics for the Jacobi polynomials via uniform treatment of Darboux's method, using not the results but the ideas discussed above.
Later in 2013,
  Liu,   Wong and  Zhao \cite{Liu-Wong-Zhao-2013}  apply the uniform Darboux's method to analyze  the Heisenberg polynomials.
  The asymptotic expansion obtained involves the Kummer function and its derivative.

We mention yet another situation, which may be termed a  singular case. We take as an  example the Pollaczek polynomials
\begin{equation}\label{Pollaczek}
 (1-ze^{i\theta})^{-\frac 1 2+ih(\theta)}(1-ze^{-i\theta})^{-\frac 1 2-ih(\theta)}=\sum^\infty_{n=0}P_n(x;a,b) z^n,
\end{equation} each of the factors in the generating
function reduces to $1$ for $z=0$, where
$h(\theta)=\frac {a\cos\theta +b}{2\sin\theta}$, $a>|b|$.
 There are two
singularities   $z=e^{\pm
i\theta}$ on the circle of convergence,
coalescing with each other when $\theta\rightarrow 0$.
 Meanwhile the exponent $h(\theta)\sim \frac {a+b} {2\theta}$ as  $\theta\rightarrow 0$, demonstrating a singular behavior. The reader is referred to \cite{Bo-Wong-1996} for an asymptotic analysis of these polynomials using integral methods.  Such types of generating function with varying exponents, regular or singular, seem  to be of interest.

\section{Difference equations}\label{sec: difference}

As remarked earlier in Section\;\ref{sec: introduction}, the development of the
asymptotic theory for difference equations took a halt in the thirties  of the last century, and did not make any progress until the turn of this century. It is interesting to note that in his lecture, given  during the conference in honor of his 80th birthday, Frank Olver expressed that
{\it In my view Birkhoff and Trjitzinsky \cite{Birkhoff-Trjitzinsky-1933} set back all research into the asymptotic solutions of difference equations for most of the 20th century.}
  Also, in 1985,  Wimp and  Zeilberger \cite{Wimp-Zeilberger-1985} made the interesting remark that
     {\it Once on the forefront of mathematical research in America, the asymptotics of the solutions of linear recurrence equations is now almost forgotten, especially by the people who need it most, namely combinatorics and computer scientists.}
Later, in 1991, Wimp \cite{Wimp-1991} also made the statement that
    {\it There are still vital matters to be resolved in asymptotic analysis. At least one
widely quoted theory, the asymptotic theory of irregular difference equations
expounded by G. D. Birkhoff and W. R. Trjitzinsky \cite{Birkhoff-1930,Birkhoff-Trjitzinsky-1933}  in the early 1930's, is
vast in scope; but there is now substantial doubt that the theory is correct in all
its particulars. The computations involved in the algebraic theory alone   are truly mind-boggling.}

The above comments were certainly the driving force behind the work carried out by Wong and Li \cite{Wong-Li-1992a,Wong-Li-1992b} in the 1990's.
In order to understand the problem better, they first studied in
\cite{Wong-Li-1992a} the simple second-order difference equation
\begin{equation*}
y(n+2)+a(n) y(n+1)+b(n) y(n)=0,
\end{equation*}
where  $a(n)$ and $b(n)$ have infinite expansions of the form
\begin{equation}\label{coefficients-DE}
 a(n)\sim\sum^\infty_{s=0} \frac {a_s}{n^s}~~~~\mbox{and}~~~~
  b(n) \sim\sum^\infty_{s=0} \frac {b_s}{n^s}
\end{equation}for large values of $n$, and $b_0\not=0$.  Then in \cite{Wong-Li-1992b}, they extended their investigation to include the more general equation
\begin{equation}\label{DE}
y(n+2)+n^p a(n) y(n+1)+n^q b(n) y(n)=0,
\end{equation}where $p$ and $q$ are integers, and
$a(n)$ and $b(n)$ are as given in
\eqref{coefficients-DE} with the leading coefficients $a_0$ and  $b_0$ being nonzero.

The basic technique adopted in \cite{Wong-Li-1992a,Wong-Li-1992b} is essentially the method of successive approximations, which is customarily used in the asymptotic theory of differential equations.  The ultimate goal of Wong's investigation was  to develop a turning point theory for the three-term recurrence relation
\begin{equation}\label{DE-x}
y_{n+1}=(a_n x+ b_n) y_n-c_n y_{n-1},~~~~n=1,2,\cdots,
\end{equation}
where $a_n$, $b_n$ and $c_n$ are constants. If $x$ is a fixed number, then the recurrence relation is equivalent to the second-order linear difference equation \eqref{DE}.  The importance of developing of a turning
point theory for the recurrence relation \eqref{DE-x} lies in the fact that many special functions of mathematical physics (Bessel functions, parabolic cylinder functions, Legendre functions, etc.) satisfy such an equation. In fact, any sequence of orthogonal polynomials satisfies an equation of the form \eqref{DE-x}; see \cite[p.\;42]{Szego-1975}.  This statement, of course, applies to classical orthogonal polynomials such as Hermite  $H_n(x)$, Laguerre $L^{(\alpha)}_n(x)$ and Jacobi $P_n^{(\alpha,\beta)}(x)$; more importantly, it includes all those that do not satisfy any second-order linear differential equations, e.g., Charlier $C_n^{(a)}(x)$,  Meixner $m_n(x; \beta, c)$, Pollaczek $P_n(x; a, b)$, Meixner-Pollaczek $M_n(x;\delta,\eta)$ and Krawtchouk polynomials $K_n^N(x; p,q)$.

While a turning-point theory for second-order differential equations has been satisfactorily developed from 1930 to 1960 by Langer \cite{Langer-1931},
Cherry \cite{Cherry-1949}, Olver \cite{Olver-1956} and others, not much progress has  been made in the development of a corresponding theory for the recurrence equation \eqref{DE-x}.  A possible explanation is that difference equations are considerably more complicated to analyze than differential
equations; see a remark made by Iserles \cite[p.\;743]{Iserles-2000}.

To illustrate the difficulty in hand, we consider the Bessel function $J_\nu(\nu x)$, which has been used as  a typical example in some of the most profound research in asymptotic analysis; see, e.g., \cite{Chester-Friedman-Ursell-1957},  \cite{Olver-1954a} and \cite{Olver-1954b}.  The differential equation approach proceeds with the
equation
\begin{equation}\label{Bessel}
\frac {d^2 w}{dx^2}=\left\{  \nu^2 \frac {1-x^2}{x^2}-\frac 1 {4x^2}\right \} w,
\end{equation}which is satisfied by $x^{\frac 1 2}J_\nu(\nu x)$.   From \eqref{Bessel}, it is  clear that $x=\pm 1$ are two turning points, i.e., zeros of the coefficients  function   multiplied  by the large parameter $\nu^2$. An application  of the transformations
\begin{equation}\label{transformation-variable}
\frac 2 3 \zeta^{\frac  3 2} =\ln \frac {1+(1-x^2)^{\frac 1 2}} x-(1-x^2)^{\frac 1 2}
\end{equation}and
\begin{equation}\label{transformation-function}
W=\hat{f}^{\frac 1 4} w,~~~~~~\hat{f}=\frac {1-x^2}{x^2\zeta}
\end{equation}
gives
\begin{equation*}
\frac {d^2 W}{d\zeta^2}=\left\{  \nu^2 \zeta+\psi(\zeta)\right\} W ,
\end{equation*}where
\begin{equation}\label{psi-def}
\psi(\zeta)=\frac 5 {16\zeta^2}+\frac {\zeta x^2(x^2+4)}{4(x^2-1)}.
\end{equation}The basic approximation equation is
\begin{equation}\label{Bessel-approx}
\frac {d^2 W}{d\zeta^2}=   \nu^2 \zeta W;
\end{equation}
two linearly independent solutions are the Airy functions
$\Ai\left(\nu^{\frac 2 3}\zeta\right )$ and $\Bi\left(\nu^{\frac 2 3}\zeta\right )$. After identifying $J_\nu(\nu x)$ and $\Ai\left(\nu^{\frac 2 3}\zeta\right )$ and matching their behavior as $x\to\infty$, one obtains the expansion
\begin{equation}\label{J-nu-nu-expan}
J_\nu(\nu x)\sim \left (\frac {4\zeta}{1-x^2}\right )^{\frac 1 4}\left\{
\frac{\Ai\left(\nu^{\frac 2 3}\zeta\right )}{\nu^{\frac 1 3}}\sum^\infty_{s=0} \frac{A_s(\zeta)}{\nu^{2s}}+
\frac{\Ai'\left(\nu^{\frac 2 3}\zeta\right )}{\nu^{\frac 5 3}}\sum^\infty_{s=0} \frac{B_s(\zeta)}{\nu^{2s}}
\right\}
\end{equation}as $\nu\to\infty$, uniformly with respect to $x\geq 0$.  The coefficients
 $A_s(\zeta)$ and $B_s(\zeta)$ are defined recursively  by $A_0(\zeta)=1$,
\begin{equation}\label{B-s}
B_s(\zeta)=\frac 1 {2\zeta^{\frac 1 2}}\int^\zeta_0 \left\{ \psi(v)  A_s(v)- A_s''(v)\right\}  \frac {dv}{v^{\frac 1 2 }}
\end{equation}and
\begin{equation}\label{A-s}
A_{s+1}(\zeta)=-\frac 1 2 B_s'(\zeta)+\frac 1 {2 }\int    \psi(\zeta)  B_s(\zeta)  d\zeta +\mbox{constants}.
\end{equation}

From a difference-equation point of view, it is natural to suggest that the same result can be obtained directly from the three-term recurrence equation
\begin{equation}\label{J-nu-difference-eq}
J_{\nu+1}(x)-\frac {2\nu} x J_{\nu}(x)+J_{\nu-1}(x)=0.
\end{equation}
This problem turns out to be considerably more difficult to tackle than what we would have anticipated.  To start with the very form of expansion \eqref{J-nu-nu-expan} immediately raises the following questions:
\begin{description}
  \item (i) Since the Airy function $\Ai(x)$ does not satisfy any second-order difference equation, how is it going to arise from the three-term recurrence relation \eqref{J-nu-difference-eq}?
  \item (ii) Since we can not change the discrete variable $n$ (or $\nu$) to a new independent variable, like what we have done with the Langer transformation \eqref{transformation-variable}, how are we going to produce the function $\zeta$ in \eqref{J-nu-nu-expan} (cf. \eqref{transformation-variable} )?
  \item (iii) Even if we can derive an asymptotic expansion similar to that in \eqref{J-nu-nu-expan}, can the coefficients be determined recursively like \eqref{B-s}-\eqref{A-s}?
  \item (iv) Can the region of validity of our expansion be as large as the one for    \eqref{J-nu-nu-expan}, i.e., the whole positive real axis $x\geq 0$; or is it just a finite interval $0\leq x\leq M$ like in the case of integral approach (see \cite[p.\;610]{Chester-Friedman-Ursell-1957})?
\end{description}
Answers to these questions are given in \cite{Wang-Wong-2002}.

In 1967, Dingle and Morgan \cite{Dingle-1967a, Dingle-1967b} made an attempt to study second-order difference equations of the form
\begin{equation}\label{Dingle-Difference-eq}
f_{n+\omega}(x)+f_{n-\omega}(x)=2\sigma(n,x) f_n(x),
\end{equation}where $\sigma(n,x)$ is a slowly varying function of the independent variable $n$. For convenience,
they expressed $\sigma(n,x)$ as
 \begin{equation}\label{Dingle-coeff-expan}
 \sigma(n,x)= a(n,x)+\omega^2 b(n,x)+\omega^4 c(n,x)+\cdots
 \end{equation}with the understanding that  $b(n,x)$ is two ``order'' smaller than $a(n,x)$, and that $c(n,x)$ is four ``order'' smaller than $a(n,x)$ and so on; see \cite[eq.\;(35)]{Dingle-1967b}.  In \eqref{Dingle-coeff-expan}, $\omega$ is used as an ``ordering'' parameter which picks out terms of comparable magnitude.  These authors called the zeros $n_0(x)$  of the equation
 \begin{equation}\label{Dingle-turing-point}
 a(n,x)=1
 \end{equation}turning points, and suggested that near $a=1$ (i.e., in a neighborhood of $n_0(x)$), two linearly independent asymptotic expansions are given by
 \begin{equation}\label{Dingle-f-1}
 f^{(1)}_n (x)=J_{n+F(x)-n_0(x)}(F(x))
 \end{equation}and
\begin{equation}\label{Dingle-f-2}
 f^{(2)}_n (x)=Y_{n+F(x)-n_0(x)}(F(x)),
 \end{equation}where
\begin{equation}\label{Dingle-F}
F(x)=\left(\frac {d a(n,x)}{dn}\right )^{-1}_{n=n_0(x)} .
\end{equation}
In the case of \eqref{J-nu-difference-eq}, $\omega=1$ and $\sigma(n, x)=n/x$. Thus, in \eqref{Dingle-coeff-expan} we just take $a(n,x)=\sigma(n, x)=n/x$ and $b(n,x)=c(n,x)=\cdots=0$.  The turning point (in the sense of Dingle and Morgan) is  at $n=n_0(x)= x$, and the function in \eqref{Dingle-F} is given by $F(x)=x$. As a result, the asymptotic solution $f^{(1)}_n (x)$ in \eqref{Dingle-f-1}
is just the function $J_{n}(x)$ itself, which is what we wish to approximate.  In any case, the arguments in \cite{Dingle-1967a,Dingle-1967b} are too sketchy and non-rigorous. Nevertheless, in view of \eqref{J-nu-nu-expan}, the result of Dingle and Morgan is most likely correct.

Another paper that is relevant to our discussion here is by Costin and Costin \cite{Costin-Costin-1996}. These authors have studied the asymptotic behavior of the solutions of recurrence relations of the form
\begin{equation}\label{Costin-difference-eq}
a_2(k\varepsilon, \varepsilon) y_{k+2}+ a_1(k\varepsilon, \varepsilon) y_{k+1}+a_0(k\varepsilon, \varepsilon) y_{k}=0
\end{equation}as $\varepsilon\to 0^+$, where the coefficients $a_i(x,\varepsilon)$, $i=0,1,2$, are $C^\infty$-functions in $x$ and $\varepsilon$ in some domain $I\times [0,  \varepsilon_0]$, $I$ being a compact interval.  Furthermore, they assume that for any $n\geq 0$,
\begin{equation}\label{Costin-coeff-expan}
a_i(x, \varepsilon)=\sum^n_{s=0} a_{i,s} (x) \varepsilon^s+O\left ( \varepsilon^{n+1}\right ),~~~~i=0,1,2,
\end{equation}
where $a_{i,s} (x)$ are $C^\infty$-functions in $x$ and  the $O$-symbol is uniform with respect to $x$. Let $\lambda_1(x)$ and $\lambda_2(x)$ be the roots of the characteristic polynomial
\begin{equation}\label{Costin-char-eq}
a_2(x, 0) \lambda^2+ a_1(x, 0) \lambda+a_0(x, 0)  =0,
\end{equation}
and assume that $\lambda_1(0)=\lambda_2(0)$
but $\lambda_1(x)\not=\lambda_2(x)$ when $x\not=0$.    Fix $\frac 1 3<\alpha<\frac 1 2$.
One of the major results in \cite{Costin-Costin-1996} is that for $k<\varepsilon^{-\alpha}$, equation
\eqref{Costin-difference-eq}  has two solutions of the form
\begin{equation*}
\exp\left\{ F_\pm\left(k\varepsilon^{\frac 1 3}, \varepsilon^{\frac 1 3}\right )\right\},
\end{equation*}
where  $\exp\left\{   F_\pm\left(x, 0\right )\right\}=\Ai(\Theta x)\pm \Bi(\Theta x)$ and $\Theta$ is an explicitly given constant in terms of $a_i(0,0)$ and $D_x a_i(0,0)$. A particular solution has the behavior
\begin{equation}\label{Costin-solution}
y_{k, \varepsilon}\sim \Ai\left(\Theta k\varepsilon^{\frac 1 3} \right)\left [ 1+
 \varepsilon^{\frac 1 3}A_1\left(\Theta k\varepsilon^{\frac 1 3}\right )+\cdots\right ]
\end{equation}for large $k<\varepsilon^{-\alpha}$. A comparison of equations \eqref{J-nu-difference-eq} and \eqref{Costin-difference-eq} readily shows that the results of Costin and Costin can not be used to derive a uniform asymptotic expansion for $J_\nu(\nu x)$ such as the one given in \eqref{J-nu-nu-expan}.

We are now ready to present a summary of the investigation carried out by Wong and his collaborators in the last twenty years.  Returning to equation \eqref{DE-x}, we first define a sequence $\{K_n\}$ recursively by
$K_{n+1}/K_{n-1}=c_n$, with $K_0$ and $K_1$ depending on the particular sequence $\{y_n\}$ satisfying \eqref{DE-x}.
Then we put $A_n=a_n K_n/K_{n+1}$,  $B_n=b_n K_n/K_{n+1}$ and $x_n=y_n/K_n$, so that \eqref{DE-x}
becomes
\begin{equation}\label{DE-x-sym}
x_{n+1}-\left (A_n x+B_n\right ) x_n+x_{n-1}=0.
\end{equation}
The coefficients $A_n$ and $B_n$ are assumed to have asymptotic expansions of the form
\begin{equation}\label{De-x-sym-coef}
A_n\sim n^{-\theta} \sum^\infty_{s=0} \frac{\alpha_s}{n^s}~~~~~~\mbox{and}~~~~~~
B_n\sim   \sum^\infty_{s=0} \frac{\beta_s}{n^s},
\end{equation}where $\theta$ is a real number and $\alpha_0\not=0$.

Let $\tau_0$ be  a constant and put $\nu:=n+\tau_0$. Clearly, the expansions in
\eqref{De-x-sym-coef} can be recast  in the form
\begin{equation}\label{De-x-recast-coef}
A_n\sim \nu^{-\theta} \sum^\infty_{s=0} \frac{\alpha'_s}{\nu^s}~~~~~~\mbox{and}~~~~~~
B_n\sim   \sum^\infty_{s=0} \frac{\beta'_s}{\nu^s}.
\end{equation}
In \eqref{DE-x-sym}, we now let $x=\nu^\theta t$ and $x_n=\lambda^n$.  Substituting \eqref{De-x-recast-coef} into \eqref{DE-x-sym} and letting $n\to\infty$ (and hence $\nu\to\infty$), we obtain the characteristic equation \begin{equation}\label{char-eq-new}
\lambda^2-(\alpha'_0t +\beta'_0)\; \lambda+1=0.
\end{equation}The roots of this equation are given by
\begin{equation}\label{roots-char-eq}
\lambda=\frac 1 2 \left [ (\alpha'_0t+\beta'_0)\pm \sqrt{(\alpha'_0t+\beta'_0)^2-4}\; \right ],
\end{equation}and they coincide when $t=t_\pm$, where
\begin{equation}\label{cr-value}
 \alpha'_0t_\pm +\beta'_0=\pm 2.
\end{equation}
The values $t_\pm$ play an important role in the asymptotic theory of the three-term recurrence relation \eqref{DE-x-sym}, and they correspond to the transition points (i.e., turning points and poles) occurring in differential equations;  \cite[p.\;362]{Olver-1974}.  For this reason, we shall also call them {\it{transition points}}. Since $t_+$ and $t_-$ have different values, we may restrict ourselves to just consider the case $t=t_+$.  In terms of the exponent $\theta$ in \eqref{De-x-recast-coef} and transition point $t_+$, we have three cases to consider; namely, (i) $\theta\not=0$ and $t_+\not=0$; (ii) $\theta\not=0$ and $t_+ =0$; and (iii) $\theta=0$. In the following, we present the main result in each of these three cases in the form of a theorem.

 {\bf Case (i)}\;  In this case, we assume $\theta>0$. The analysis for the case $\theta<0$ is very similar; for an important example with $\theta=-1$, see \cite{Wang-Wong-2002}. Since the discussion of the two cases $t_+>0$ and $t_+<0$ are similar, we shall consider only the case $t_+>0$. Furthermore, to make the presentation simpler, we also assume that $t_-<0$. These two assumptions together are equivalent to the condition $|\beta_0|<2$.  Before stating our first theorem, we need to impose one more condition; that is,
\begin{equation}\label{alpha-1-beta-1}
 \alpha_1'=\beta_1'=0.
\end{equation}
It is interesting to note that this condition holds in most of the classical cases.  In fact, in \cite{Dingle-1967a} Dingle and Morgan have assumed that $\alpha_{2s+1}=\beta_{2s+1}=0$ for $s=0,1,2,\cdots.$  For a result without this assumption, we refer the reader to \cite{Wang-Wong-2003}. Now we define the function

\begin{equation}\label{conformal-mapp}
\left\{\begin{array}{ll}
\displaystyle{\frac 2 3 \left [ \zeta(t)\right ]^\frac 3 2 =\alpha'_0 t^{\frac 1 \theta} \int_{t_+}^t \frac {s^{-1/\theta} ds}
{\sqrt{(\alpha'_0 s+\beta'_0)^2-4}}-\ln \frac { \alpha'_0 t+\beta'_0 + \sqrt{(\alpha'_0t+\beta'_0)^2-4}} 2,} & t\geq t_+, \\[.5cm]
\displaystyle{\frac 2 3 \left [-\zeta(t)\right ]^\frac 3 2 =-\alpha'_0 t^{\frac 1 \theta} \int_t^{t_+} \frac {s^{-1/\theta} ds}
{\sqrt{4-(\alpha'_0 s+\beta'_0)^2}}+\arccos  \left( \frac { \alpha'_0 t+\beta'_0  } 2\right ) }, & t< t_+,
       \end{array}
 \right .
\end{equation}
which plays the role of the Langer transform for differential equations; cf. \eqref{transformation-variable}.

\begin{thm}
Assume that the coefficients $A_n$ and $B_n$ in the recurrence relation \eqref{DE-x-sym} have asymptotic expansions of the form given in \eqref{De-x-sym-coef} with $\theta\not=0$ and $|\beta_0|<2$. (See the statement just before \eqref{alpha-1-beta-1}.) Let $\zeta(t)$ be defined as in \eqref{conformal-mapp},  and recall the number $\nu:=n+\tau_0$ in the asymptotic expansions in \eqref{De-x-recast-coef}.  Then, with $x=\nu^\theta t$, equation \eqref{DE-x-sym} has a pair of solutions $P_n(x)$ and $Q_n(x)$ given by
\begin{equation*}
P_n\left (\nu^\theta t\right )\sim \left (\frac {4\zeta}{(\alpha'_0 t+\beta'_0)^2-4 }\right )^{\frac 1 2}
\left [
 \Ai\left(\nu^{\frac 2 3}\zeta\right ) \sum^\infty_{s=0} \frac{\widetilde A_s(\zeta)}{\nu^{s-\frac 1 6}}+
 \Ai'\left(\nu^{\frac 2 3}\zeta\right ) \sum^\infty_{s=0} \frac{\widetilde B_s(\zeta)}{\nu^{s+\frac 1 6}}
\right  ]
\end{equation*}
and
\begin{equation*}
Q_n\left (\nu^\theta t\right )\sim \left (\frac {4\zeta}{(\alpha'_0 t+\beta'_0)^2-4 }\right )^{\frac 1 2}
\left [
 \Bi\left(\nu^{\frac 2 3}\zeta\right ) \sum^\infty_{s=0} \frac{\widetilde A_s(\zeta)}{\nu^{s-\frac 1 6}}+
 \Bi'\left(\nu^{\frac 2 3}\zeta\right ) \sum^\infty_{s=0} \frac{\widetilde B_s(\zeta)}{\nu^{s+\frac 1 6}}
\right  ],
\end{equation*}
where the coefficients $\widetilde A_s(\zeta)$ and $\widetilde B_s(\zeta)$ are determined successively from some recursive formulas, beginning with  $\widetilde A_0(\zeta)=1$ and $
\widetilde B_0(\zeta)=0$.
\end{thm}

{\bf Case (ii)}\;  Historically, this case was dealt with nearly 10 years after Case (iii) was settled. One of the reasons is that it was not clear what exactly the dominant approximant should be, although the example of Laguerre polynomials suggested that it most likely is a Bessel function. Since $\theta\not=0$ in this case, we may let $\tau_0:=-\alpha_1/(\alpha_0\theta)$ and $N:=n+\tau_0$.  Then, from \eqref{DE-x-sym} and \eqref{De-x-sym-coef}, we have
\begin{equation}\label{coefficients-(ii)-DE}
A_nx+B_n\sim n^{-\theta} \sum^\infty_{s=0} \frac{\alpha_s}{n^s} x +   \sum^\infty_{s=0} \frac{\beta_s}{n^s}
:=\sum^\infty_{s=0} \frac{\alpha'_s t +\beta'_s}{N^s},
\end{equation}
where $x:=N^\theta t$.  A simple calculation gives
\begin{equation}\label{coefficients-leading-(ii)-DE}
\alpha'_0=\alpha_0,~~~~\alpha'_1=0,~~~~\beta'_0=\beta_0,~~~~\beta'_1=\beta_1,~~~~\beta'_2=\beta_2+\beta_1\tau_0.
\end{equation}
If $\beta_0=2$ (or $-2$), from \eqref{cr-value} it follows that one of the transition points is zero. Without loss of generality, we assume that $t_+=0<t_-$ and $\beta_0=2$, $\alpha_0<0$. (For other cases, see Remark \ref{rem:case-beta=-2} below.) As in Case (i), we assume that
$\beta'_1=\beta_1=0$ so that
\begin{equation}\label{cr-t+}
\alpha'_1t_+ +\beta'_1=0.
\end{equation}
(In most of the classical cases, $\beta'_1=\beta_1=0$; see the statement following \eqref{alpha-1-beta-1}.)

With $\alpha'_0<0$ and $\beta_0=2$,
we now define the function
\begin{equation}\label{conformal-map-ii}
\pm \zeta(t)= \arccos \left (  \frac { \alpha'_0 t+\beta'_0  } 2 \right )+ \alpha'_0 t^{\frac 1 \theta}   \int_a^{t} \frac {\phi^{-1/\theta}}
{\sqrt{4-(\alpha'_0 \phi+\beta'_0)^2}} d\phi
\end{equation}
for $ t< t_-$, where ``$+$'' sign is taken for $\theta<2$ and ``$-$'' sign is taken for $\theta>2$. The sign is chosen so that $\zeta(t)$ is positive for $t\in (0, t_-)$.  The choice of the lower limit of integration is just for the purpose of convergence of the integral. For instance, we can choose $a=0$ if $\theta <0$ or $\theta>2$ and
\begin{equation}\label{a-def-ii}
a=\left\{
\begin{array}{ll}
  -\infty,  & t<0, \\[.2cm]
  t_-, & 0\leq t\leq t_-
\end{array}\right . ~~~~~~\mbox{if}~0<\theta<2.
\end{equation}
With this choice, it can be shown that $\zeta(0)=0$.

\begin{thm}\label{thm:De-case-ii}
Assume that the coefficients $A_n$ and $B_n$ in the recurrence relation \eqref{DE-x-sym} have asymptotic expansions given in \eqref{De-x-sym-coef} with $\theta\not=0, 2$ and $\beta_0=2$.  Let $t_+=0$ be a transition point defined in \eqref{cr-value}, and the function $\zeta(t)$ be given as in \eqref{conformal-map-ii} and \eqref{a-def-ii}. Then,   \eqref{DE-x-sym} has a pair of linearly independent solutions
\begin{equation}\label{P-n-ii}
P_n\left (N^\theta t\right )\sim N^{\frac 1 2}\left (\frac {4\zeta^2}{4-(\alpha'_0 t+\beta'_0)^2 }\right )^{\frac 1 4}
\left [
 J_\nu\left(N\zeta\right ) \sum^\infty_{s=0} \frac{\widetilde A_s(\zeta)}{N^{s}}+
 J_{\nu+1}\left(N\zeta\right ) \sum^\infty_{s=0} \frac{\widetilde B_s(\zeta)}{N^{s}}
\right  ]
\end{equation}
and
\begin{equation*}
Q_n\left (N^\theta t\right )\sim N^{\frac 1 2}\left (\frac {4\zeta^2}{4-(\alpha'_0 t+\beta'_0)^2 }\right )^{\frac 1 4}
\left [
 W_\nu\left(N\zeta\right ) \sum^\infty_{s=0} \frac{\widetilde A_s(\zeta)}{N^{s}}+
 W_{\nu+1}\left(N\zeta\right ) \sum^\infty_{s=0} \frac{\widetilde B_s(\zeta)}{N^{s}}
\right  ]
\end{equation*}for $-\infty<t\leq t_--\delta$, $\delta$ being an arbitrary positive constant. Here,
 $W_\nu(x):=Y_\nu(x)-i J_\nu(x)$, $x=N^\theta t$, $N=n+\tau_0$, $\tau_0=-\alpha_1/(\alpha_0\theta)$ and $\nu$ is given by
 \begin{equation}\label{nu-def-ii}
 \nu=\sqrt{\frac{1+4\beta'_2}{(\theta-2)^2}} .
 \end{equation}
 The coefficients  $\widetilde A_s(\zeta)$ and $\widetilde B_s(\zeta)$ are determined successively from some recursive formulas, beginning with  $\widetilde A_0(\zeta)=1$ and $\widetilde B_0(\zeta)=0$.
\end{thm}

In the above theorem, we have used the functions $J_\nu(x)$ and $W_\nu(x):=Y_\nu(x)-i J_\nu(x)$ as two linearly independent solutions of the Bessel equation. This is because $J_\nu(x)$ and $H^{(1)}_\nu(x)$ form a numerically satisfactory pair when $0<\arg x<\pi$; see \cite[\S 10.2(iii)]{nist-2010}. Note that $i\zeta(t)$ is negative for $t<0$; thus $J_\nu(N\zeta)$ is exponentially large and $H^{(1)}_\nu(N\zeta)$  is exponentially small as $N\to\infty$.

\begin{rem}\label{rem:case-beta=-2}
For the case $\alpha_0>0$ and $\beta_0=-2$, the two transition points $t_\pm$ satisfy $t_-=0<t_+$.  Note that we again have one transition point at the origin and the  other being positive.
Put $\mathcal{P}_n(x):=(-1)^n P_n(-x)$. Theorem \ref{thm:De-case-ii} then applies to $\mathcal{P}_n(x)$.  If the two transition points satisfy $t_-<0=t_+$, we set $x:=-N^\theta t$, instead of $x:=N^\theta t$.
Theorem \ref{thm:De-case-ii} again holds with $P_n\left (N^\theta t\right )$ replaced by $P_n\left (-N^\theta t\right )$ (in the case $\alpha_0>0$,  $\beta_0=2$) or  $(-1)^n P_n\left (-N^\theta t\right )$ (in the case $\alpha_0<0$,  $\beta_0=-2$).
\end{rem}

\begin{rem}\label{rem:case-exceptional}
In Theorem \ref{thm:De-case-ii}, we assume $\theta\not=0,2$.  The situation with $\theta=0$ was studied in \cite{Wang-Wong-2005}, and will be discussed in Case (iii) below. The result in that case holds regardless whether the transition point is at 0 or not.  $\theta=2$ is exceptional  in our present case (i.e., Case (ii)), since $\theta-2$ appears in the denominator of \eqref{nu-def-ii}.  However, $\theta=2$ and one of the transition point being zero is exactly the case of Wilson polynomials.  For the results in that case, we refer to \cite{Li-Wang-Wong-2020}.
\end{rem}

The presentation of Case (ii) here is based on the work of Cao and Li \cite{Cao-Li-2014}.
\vskip 1cm

{\bf Case (iii)}\;  We now consider the case   $\theta=0$. Here, there is no need to make the scale change
 $x=\nu^\theta t$ in Case (i) or $x=N^\theta t$ in Case (ii); cf. \eqref{De-x-recast-coef} and \eqref{coefficients-(ii)-DE}. Also, we will not be concerned with whether the transition point $t_+$($=x_+$)
in \eqref{cr-value} is at the origin or not.  However, we still need to make the assumptions
$\alpha_1=\beta_1=0$; see \eqref{alpha-1-beta-1} and \eqref{cr-t+}. Let
  \begin{equation}\label{tau-0-iii}
  \tau_0:=-\frac{\alpha_3 x_++\beta_3}{2(\alpha_2x_++\beta_2)},~~~~~~N:=n+\tau_0,\end{equation}
and define
\begin{equation}\label{nu-def-iii}
\nu=\left(\alpha'_2x_+ +\beta'_2+\frac 1 4\right ) ^{1/2},
\end{equation}
\begin{equation}\label{conformal-map-iii}
 \zeta^{\frac 1 2}(x)= \arccosh \left (  \frac { \alpha_0 x+\beta_0  } 2 \right ).
\end{equation}
Under the assumption in \eqref{alpha-1-beta-1}, it is easily verified that
$\alpha'_2=\alpha_2$ and $\beta'_2=\beta_2$; see \eqref{coefficients-leading-(ii)-DE}.

\begin{thm}\label{thm:De-case-iii}
Assume that the coefficients $A_n$ and $B_n$ in the recurrence relation \eqref{DE-x-sym} have asymptotic expansions given in \eqref{De-x-sym-coef} with $\theta=0$.    Let $x_\pm$($=t_\pm$) be the transition points defined by \eqref{cr-value}. Then,  equation \eqref{DE-x-sym} has a pair of linearly  independent solutions
\begin{equation}\label{P-n-iii}
\begin{array}{rl}
\displaystyle{ P_n\left (x\right )\sim }
& \displaystyle{\left (\frac {4\zeta}{(\alpha_0 x+\beta_0)^2 -4 }\right )^{\frac 1 4}\left [ N^{\frac 1 2}
 I_\nu\left(N\zeta^{\frac 1 2}\right ) \sum^\infty_{s=0} \frac{  A_s(\zeta)}{N^{s}}   \right .}\\[.5cm]
 &~~~~~~~~~~~~~~~~~~~~~~~~~~~~\displaystyle{\left . +
N^{\frac 1 2}\zeta^{\frac 1 2} I_{\nu-1}\left(N\zeta^{\frac 1 2}\right ) \sum^\infty_{s=0} \frac{  B_s(\zeta)}{N^{s}}
\right  ]}
\end{array}
\end{equation}
and
\begin{equation}\label{Q-n-iii}
\begin{array}{rl}
\displaystyle{
Q_n(x)\sim  }&
\displaystyle{\left (\frac {4\zeta}{(\alpha_0 x+\beta_0)^2 -4 }\right )^{\frac 1 4}\left [N^{\frac 1 2}
 K_\nu\left(N\zeta^{\frac 1 2}\right ) \sum^\infty_{s=0} \frac{  A_s(\zeta)}{N^{s}}  \right .}\\[.5cm]
 &~~~~~~~~~~~~~~~~~~~~~~~~~~~~\displaystyle{-\left .
N^{\frac 1 2}\zeta^{\frac 1 2} K_{\nu-1}\left(N\zeta^{\frac 1 2}\right ) \sum^\infty_{s=0} \frac{  B_s(\zeta)}{N^{s}}
\right  ]}\end{array}
\end{equation}as $n \to\infty$, for $x_-+\delta \leq x <\infty$, $\delta$ being an arbitrary positive constant,   where   $N=n+\tau_0$, $\tau_0$ is given in
\eqref{tau-0-iii}, $\nu$ is given in \eqref{nu-def-iii}
and    $\zeta(x)$ is defined by \eqref{conformal-map-iii}. The coefficients
  $A_s(\zeta)$
  and $B_s(\zeta)$ can be determined successively for any given  $A_0(\zeta)$ and $B_0(\zeta)$.
\end{thm}

\section{Riemann-Hilbert approach}\label{sec: Riemann-Hilbert}

A significant connection between orthogonal
polynomials and matrix-valued Riemann-Hilbert problems was established  in the work of Fokas, Its and Kitaev \cite{Fokas-Its-Kitaev}.
Assuming that there is a sequence of polynomials $\{p_n\}$, orthonormal with respect to the
weight function $d\alpha(\zeta)=w(\zeta) d\zeta$   supported  on a curve $\Gamma$, the formulation of Fokas, Its and Kitaev is the following
Riemann-Hilbert problem for a $2\times 2$ matrix-valued function $Y(z)$:
\begin{description}
    \item($Y_1$) $Y (z)$ is analytic in $\mathbb{C}\setminus \Gamma$.
  \item ($Y_2$) The jump condition on $\Gamma$ is
  \begin{equation}\label{RHP-jump}
  Y_+(\zeta)=Y_-(\zeta)\begin{pmatrix}
                         1 & w(\zeta) \\
                         0 & 1 \\
                       \end{pmatrix},~~~~\zeta\in \Gamma,
  \end{equation} where $Y_\pm(\zeta)$ denote the limits of $Y(z)$ as $z\to\zeta\in \Gamma$ from the left and right of the oriented curve $\Gamma$, respectively.
  \item ($Y_3$) The behavior at infinity is
\begin{equation}\label{RHP-infty}
Y (z)=\left (I +O\left (\frac 1 z\right )\right )
 \begin{pmatrix}
  z^n & 0 \\
0 & z^{-n} \\
  \end{pmatrix},~~~~z\to\infty .
\end{equation}
\end{description}
The significance lies in that the unique solution  to
the above Riemann-Hilbert problem
can be determined explicitly in terms of
the orthogonal polynomials, that is,
\begin{equation}\label{RHP-solution}
Y (z)=\begin{pmatrix}
        \pi_n(z) &\displaystyle{ \frac 1 {2\pi i}\int_\Gamma\frac {\pi_n(\zeta) w(\zeta) d \zeta}{\zeta-z} }\\[.35cm]
   -2\pi i\gamma_{n-1}^2 \pi_{n-1}(z) &\displaystyle{ -\gamma_{n-1}^2  \int_\Gamma\frac {\pi_{n-1}(\zeta) w(\zeta) d \zeta}{\zeta-z} }\\
      \end{pmatrix},~~~~z\in \mathbb{C}\setminus\Gamma,
\end{equation}where $\pi_n(z)$ is the monic polynomial corresponding to $p_n(z)$, and $\gamma_n$ is the leading coefficient of $p_n(z)$ such that
$p_n(z)=\gamma_n\pi_n(z)$.
The uniqueness of the solution $Y(z)$ can be justified by using Liouville's theorem, with asymptotic behavior  appropriately preassigned to each critical point   such as an endpoint  or a possible point  of intersection  of $\Gamma$, etc.

To extract  asymptotics  from \eqref{RHP-solution},
we may use   a steepest descent type
method for oscillatory
  Riemann-Hilbert problems  introduced by Deift and Zhou \cite{Deift-Zhou-1993}; see Deift et al. \cite{Deift-Kriecherbauer-McLaughlin-Venakides-Zhou-2001} for a brief historical account of how this steepest descent method, also termed the Deift-Zhou method or Riemann-Hilbert approach, was further developed.
This  new approach to study asymptotic questions for orthogonal polynomials
is rigorous, essentially global, and  has proved to be very powerful.
The pioneering works in which  the nonlinear steepest descent method was actually applied to orthogonal polynomials are
by   Deift et al. \cite{Deift-Kriecherbauer-McLaughlin-Venakides-Zhou-1999a} on orthogonal polynomials with varying exponential weight, with an application to a rigorous proof of the universality in random matrix theory,  and  by Bleher  and  Its
\cite{Bleher-Its-1999} on semi-classical asymptotics  with applications in random matrix theory; see also \cite{Deift-Kriecherbauer-McLaughlin-Venakides-Zhou-1999b}.

There are many more publications afterwards. Here we mention only a very few. In \cite{Kuijlaars-McLaughlin-Van Assche-Vanlessen-2004}, Kuijlaars et al.
 studied the polynomials orthogonal with respect to
the modified Jacobi weight. This provides  a situation where one seems to have no choice but to use the Riemann-Hilbert approach, since the orthogonal measure is the only property  known.  Discrete orthogonal polynomials are investigated in Johansson \cite{Johansson-2001} and Baik et al.
\cite{Baik-Kriecherbauer-McLaughlin-Miller-2007}. The orthogonal polynomials with respect to Freud weights   are considered in
Kriecherbauer and McLaughlin \cite{Kriecherbauer-McLaughlin-1999}.
Wong and his coauthors have made attempts to bring in
a global treatment to the Riemann-Hilbert approach; see \cite{Wang-Wong-2005b}
and a brief account in \cite{Wong-Zhao-2009}, in which
all zeros of the polynomials belong to one
single domain of uniformity.

 There have been excellent introductions to the Riemann-Hilbert approach. The reader is
referred to the books of Deift \cite{Deift-1999}  and  Fokas et al. \cite{Fokas-Its-Kapaev-Novokshenov-2006}, and   the paper of Kuijlaars
 \cite{Kuijlaars-2003},
 for a better understanding
of this method; see  also \cite{Wong-Zhao-2009} and \cite{Zhao-2006} for some more recent progress.

The asymptotic analysis of the Riemann-Hilbert problem consists of the following  successive transformations, starting from the Riemann-Hilbert problem for $Y$.
\begin{description}
  \item (i) The first transformation, if applicable,  is a re-scaling of the variable to make the equilibrium measure compactly supported.
  \item (ii) The second transformation is a normalization of the behavior at infinity,  involving  a function $g$ which is the  log transformation of the equilibrium measure, and turning the problem into  an oscillatory Riemann-Hilbert problem.
  \item (iii) The next transformation is the central piece of the analysis, which
 involves a factorization of the jump matrix and a deformation of the contour; oscillatory terms are transformed into exponentially decaying terms and  may be
neglected along the deformed contour, which reminds us of Debye's
classical method of steepest descent for integrals.
\item (iv) The last transformation combines together  the parametrices at all critical points such as the endpoints of the contour, and edges of the support of the equilibrium measure. In many circumstances
 technical difficulties may arise in the construction  of such a   parametrix, that is, an explicitly given local solution to the  modified  Riemann-Hilbert problem in a neighborhood of the critical point.
\end{description}

Now we
  distinguish two types of  parametrices.
For the first type, an example is the  frequently used parametrix
involving
the Airy functions (\cite{Deift-1999, Forrester-1993, Tracy-Widom-1994a}),   usually constructed in a neighborhood of the soft edge of the spectrum. 
 An application of it to the random matrix theory is to describe the
soft edge behavior by using  the    Airy kernel
\begin{equation*}
 \mathbb{A}(x, y) = \frac {\Ai(x)\Ai'(y)-\Ai(y)\Ai'(x)}{x-y}.
\end{equation*}
While the Bessel kernel,
\begin{equation*}
\mathbb{J}_\alpha(x,y)=\frac{J_\alpha(\sqrt x)\sqrt y J'_\alpha(\sqrt y)-
J_\alpha(\sqrt y)\sqrt x J'_\alpha(\sqrt x)}{2(x-y)},
\end{equation*}
resulted  from a Bessel parametrix, is often used to  characterize
  the behavior
at the hard edge of the spectrum; cf. \cite{Forrester-1993, Kuijlaars-McLaughlin-Van Assche-Vanlessen-2004, Tracy-Widom-1994b}.
The list of these parametrices is much longer, including  special functions such as the confluent hypergeometric functions, the parabolic cylinder functions, the Painlev\'{e} functions, and more general Painlev\'{e}-type functions.  There is a quite vast literature on these
 more transcendental
  parametrices; see, e.g., \cite{Fokas-Its-Kapaev-Novokshenov-2006}, \cite{Xu-Dai-Zhao-2020} and \cite{Xu-Zhao-2020} and the references therein.

A common feature of the above   parametrices, is that the special function involved satisfies a certain differential equation, linear or nonlinear. The reason for this  phenomenon  is that each of these parametrices  can be transformed, by 
making use of a certain local conformal  mapping, into a model Riemann-Hilbert problem, say  $P(\zeta)$, with constant jump.

In the constant jump case, $P(\zeta)$ and $P'(\zeta)$ share the same jump condition.
Hence
all jumps vanish for
$A(\zeta):=P'(\zeta)\{P(\zeta)\}^{-1}$. Then $A(\zeta)$ is analytic in $\mathbb{C}$, with only isolated singularities of known types. This would give  a differential equation $P'(\zeta)=A(\zeta) P(\zeta) $.

However, there are still paramatrices  of another type, in which
  the model problems could not
be turned into the ones with constant jumps. In what follows,
we briefly describe three   examples of this type. These examples deal with  polynomials  orthogonal with respect to a   logarithmic weight,
a Bessel-function weight, and the Freud weight, separately.

\begin{exa} {\bf{Logarithmic weight.}}
In 2018, Conway and  Deift \cite{Conway-Deift-2018} derived the asymptotics of the recurrence coefficients
for the polynomials orthogonal with respect to the logarithmic weight
$w(x) dx$,  where
\begin{equation}\label{log-weight}
w(x)=\ln\frac {2k}{1-x},~~~~x\in [-1,1],
\end{equation}and $k>1$ so that  the measure is strictly positive on its support. \end{exa}
The formulation in \eqref{RHP-jump} and \eqref{RHP-infty} fits well with $\Gamma=[-1, 1]$, the weight $w$ specified in \eqref{log-weight}, and $Y(z)$ having at most weak singularities at the endpoints $z=\pm 1$. Following a now standard series of transformations,  we  change the original
Riemann-Hilbert problem for $Y$ into a modified  problem for  $Q$, of which the jumps
on curves adjacent to   $z=1$ are
\begin{equation}\label{log-jump-modified}
J_Q(s)=\left\{
\begin{array}{ll}
  \begin{pmatrix}
   1 & 0 \\
 \frac{F^2} w \phi^{-2n} & 1 \\
  \end{pmatrix},
  & s\in \Sigma_1\cup\Sigma_2, \\[.4cm]
 \begin{pmatrix}
   0 & 1 \\
   -1 & 0 \\
 \end{pmatrix},
  & s\in (-1, 1),\\[.4cm]
 \begin{pmatrix}
   1 & 0 \\
 \left(\frac{F^2} {w_+}+\frac{F^2} {w_-} \right ) \phi^{-2n} & 1 \\
  \end{pmatrix}, & s\in (1, +\infty);
\end{array}\right .
\end{equation}cf. \cite{Conway-Deift-2018}.
Here,  $(~)_\pm$ denote the boundary values (cf. \eqref{RHP-jump}), all curves are in a rightward direction,  $\Sigma_1$ and $\Sigma_2$ are respectively the upper and lower boundaries of the lens opening from $(-1, 1)$,
 $w(z)=\ln\frac {2k}{1-z}$ is analytic in $\mathbb{C}\setminus [1, +\infty)$,   the mapping $\phi(z)$ and the Szeg\H{o} function $F(z)$ (cf. \cite{Kuijlaars-McLaughlin-Van Assche-Vanlessen-2004})  defined by
 \begin{equation*}
 \phi(z)=z+\sqrt{z^2-1},~~~~F(z)=\exp\left(
 \frac {(z^2-1)^{1/2}}{2\pi i}\int^1_{-1}\frac{\ln w(s)}{(s^2-1)_+^{1/2}}\frac {ds}{s-z}\right )
 \end{equation*}are analytic in $\mathbb{C}\setminus [-1, 1]$.
The new difficulty   is that no explicit local solution is known
near the logarithmic singularity $x=1$ of the weight \eqref{log-weight}, since the jump matrices in \eqref{log-jump-modified} are   unlikely to be simultaneously turned
into constant ones.

Fortunately,   for the above problem
with logarithmic weight, no local solution is needed.
The strategy in \cite{Conway-Deift-2018} is to compare the above modified  Riemann-Hilbert problem with the corresponding
modified Riemann-Hilbert problem for the Legendre case. The comparison is not to  examine the quotient, but to estimate the difference   of solutions to
 these two Riemann-Hilbert problems.
The estimates given in \cite{Conway-Deift-2018} seem  to allow  an  effective
comparison of two Riemann-Hilbert problems on the same contour in the general case.
\vskip .3cm

\begin{exa}
{\bf{Bessel function weight}}.  In 2016,   Dea\~{n}o, Kuijlaars  and  Rom\'{a}n
\cite{Deano-Kuijlaars-Roman-2016}
investigated the asymptotics of the polynomials $P_n(z)$  that are orthogonal with respect
to the weight function $J_\nu(x)$  on $[0, \infty)$,  where $J_\nu$ is the Bessel function of order $\nu\geq 0$.
The Bessel function is oscillatory with an amplitude that decays like $O(x^{-1/2})$ as
$x\to+\infty$, and therefore the moments
$\int_0^\infty x^j J_\nu(x) dx$
do not exist. However, the polynomials $P_n$ can be defined  via a regularization of the
weight with an exponential factor; see
Asheim and  Huybrechs \cite{Asheim-Huybrechs-2013}.
\begin{equation}\label{Bessel-OP-def}
\int_0^\infty P_n(x; s) x^j J_v(x) e^{-sx} dx=0,~~~~j=0,1,\cdots, n-1,
\end{equation}and then
\begin{equation}\label{Bessel-OP-lim}
P_n(x)=\lim_{s\to 0^+} P_n(x; s).
\end{equation}
\end{exa}
From   numerical experiments, Asheim
and Huybrechs \cite{Asheim-Huybrechs-2013} observe  that the zeros of $P_n(z)$ seem to cluster along the vertical line
$\Re z =\nu\pi/2$. When $0\leq \nu\leq 1/2$, Dea\~{n}o, Kuijlaars  and  Rom\'{a}n
\cite{Deano-Kuijlaars-Roman-2016} derive  large-$n$ asymptotic behavior
 of $P_n(in\pi z)$ by using the Riemann-Hilbert approach, and  give a rigorous proof of the observation in this case.
 More precisely, they have shown that  for
the scaled zeros, that is, dividing the imaginary parts of the zeros by $n$ while keeping the real
parts fixed,  the limiting curve is a vertical line segment of $\Re z =\nu\pi/2$.
They further conjecture, based also on numerical evidence, that
  there is a limiting curve for
the scaled zeros that differs from the vertical line segment  when  $\nu>1/2$,
   partially because
   in this case the method in \cite{Deano-Kuijlaars-Roman-2016} to construct
a local parametrix at the origin fails: this difficulty may very well be related to the
different behavior of the zeros.
Yet
a very recent research \cite{Zhao-Wu-Zhao-2021} shows that for $\nu>1/2$, the limiting curve for
the scaled zeros is still a vertical   segment  of  $\Re z =\nu\pi/2$.

The most technical part of  \cite[(3.50)-(3.51),\;Prop.\;3.16]{Deano-Kuijlaars-Roman-2016} is the analysis of the local parametrix at the origin, with jump matrices involving   the Bessel functions in a complicated way. Since by no means one can transform the jumps into   constant matrices, an integral operator argument was brought in to show that the modified local Riemann-Hilbert problem has a solution for $n$ large enough.
The leading order
behavior of the parametrix has not yet been explicitly given, though it does not affect  the determination  of zero distribution.

In 1982, Wong \cite{Wong-1982}  considered quadrature formulas for oscillatory integral transforms with Fourier kernel and Bessel kernels.
Also,  in a paper \cite{Dai-Hu-Wang-2015}, Dai, Hu and Wang have proposed a study of the  asymptotics of orthogonal polynomials whose
weight function  is closely related to the Bessel weight stated above.
\vskip .3cm

\begin{exa}\label{exa:freud-weight}
{\bf{Freud weight}}. In 1999,
 Kriecherbauer and McLaughlin \cite{Kriecherbauer-McLaughlin-1999} used the Riemann-Hilbert approach to study the strong asymptotics of the polynomials orthogonal with respect to the
 Freud weight
\begin{equation}\label{freud-weight}w_\beta(x) dx =e^{-\kappa_\beta |x|^\beta} dx,~~x\in \mathbb{R},~~\kappa_\beta =  \frac {\Gamma(\frac \beta 2)\Gamma(\frac 1 2)} {\Gamma(\frac {\beta+1} 2)}, ~~\beta>0.\end{equation}
Although the Riemann-Hilbert approach works, an obstacle arises, that is, to find  the explicit local solution  in a  neighborhood of the origin.
\end{exa}
When $0<\beta <1$, the problem is reduced to solving
a model  Riemann-Hilbert problem 
for a certain $2\times 2$ matrix-valued function $L$:
\begin{align}
   &  L: \mathbb{C}\backslash \mathbb{R}\rightarrow \mathbb{C}^{2\times 2}~ \mbox{is~analytic}, \nonumber\\
  & L_+(s)=L_-(s) v_L(s)~\mbox{for} ~s\in \mathbb{R},\label{L-jump}\\
  & L(s)=O(\ln |s|)~\mbox{as}~s\rightarrow 0,\label{L-behavior-origin}\\
  & L(s)\rightarrow I~\mbox{for}~s\rightarrow \infty,\label{L-behavior-infty}
\end{align}
where the jump matrix is given by
\begin{equation*}
v_L(s)= \left(
                   \begin{array}{cc}
                   1  & -\eta_L(s) \\
                   \eta_L(-s)   & 1\\
                   \end{array}
                 \right),~~s\in  \mathbb{R},
\end{equation*}
with
\begin{equation*}
\eta_L(s)=2i (-1)^{n+1} e^{-|s|^\beta \sin(\pi
\beta/2)}\sin\left(|s|^\beta\cos\frac{\pi\beta}2\right )
{\bf{1}}_{[0, \infty)}(s), ~s\in \mathbb{R},
\end{equation*}
${\bf{1}}_{[0, \infty)}(s)$ denoting
the indicator function   of the set $[0, \infty)$,
and $n$ being the polynomial degree; see  \cite[(6.41)-(6.43)]{Kriecherbauer-McLaughlin-1999}.

In \cite{Kriecherbauer-McLaughlin-1999}, the existence and uniqueness of $L(s)$ were established in the Wiener class.
However, as commented in Deift et al. \cite[p.\;62, Rmk.\;3.6]{Deift-Kriecherbauer-McLaughlin-Venakides-Zhou-2001},

{\it For $0<\beta<1$ the leading order behavior of the solution to the model problem at the origin has not been
determined explicitly. It is defined
through a Riemann-Hilbert problem which only depends on $\beta$ and on the parity of $n$.}

The reader is  referred to \cite{Kriecherbauer-McLaughlin-1999} and \cite{Wong-Zhao-2016} and the references there  for polynomials orthogonal
with respect to weights similar to   the Freud weight.
For example,
 Chen and Ismail \cite{Chen-Ismail-1998} considered the Freud-like orthogonal
polynomials arising from a recurrence relation related to the indeterminate
moment problems; see also Dai, Ismail  and Wang \cite{Dai-Ismail-Wang-2014}.
In fact, the polynomials orthogonal with respect to the Bessel weight in \cite{Deano-Kuijlaars-Roman-2016}; cf. \eqref{Bessel-OP-def} and \eqref{Bessel-OP-lim}, provides another example of Freud-type weights.
Expressing $J_\nu(x)$ as a linear combination of $K_\nu(\pm ix)$, one  encounters a varying exponential weight of
Freud-type $e^{-n\pi |x|}$  with a potential function $\pi|x|$.
The jump  \eqref{L-jump} with constraint \eqref{L-behavior-infty}, can not be transformed into a constant jump.
Further,
to solve the model problem, one may either  represent them in terms of  known special functions, or  alternatively  use new special functions to serve the same purpose.
\vskip .5cm

It is well known that the Riemann-Hilbert problems are closely related to singular integral equations; see, e.g., Muskhelishvili \cite{Muskhelishvili-1953}.
Assume
\begin{equation}\label{Y-to-Phi}
Y(z)=\frac 1 {2\pi i} \int_\Gamma \frac {\Phi(\zeta) d\zeta}{\zeta-z}, ~~~~z\not\in \Gamma,
\end{equation}where $\Phi(\zeta)$ for $\zeta\in \Gamma$  is the new matrix-valued unknown function. From the Plemelj formula we have
\begin{equation*}
Y_\pm(\zeta)=\pm\frac  12  \Phi(\zeta)+  \frac 1 {2\pi i}
\int_\Gamma{\!\!\!\!\!\!\!\!  {}-{} \,\,\,\, }
\frac {\Phi(\tau) d\tau}{\tau-\zeta}, ~~~~\zeta \in \Gamma,
\end{equation*}where the  bar indicates that the integral is a Cauchy principal value; cf. \eqref{Hilbert-one-side}.
Substituting the above formulas into the jump condition \eqref{RHP-jump}, we obtain a homogeneous singular integral equation
\begin{equation}\label{singular-integral-eq}
\Phi(\zeta)-\frac 1 {2\pi i}
\int_\Gamma{\!\!\!\!\!\!\!\!  {}-{} \,\,\,\, }
\frac {\Phi(\tau) d\tau}{\tau-\zeta}\begin{pmatrix}
                                      0 & w(\zeta) \\
                                      0 & 0 \\
                                    \end{pmatrix}=0
\end{equation}with index $0$.
Finding asymptotic results from \eqref{singular-integral-eq} for specific $\Gamma$  or more general integral equations is an important topic worth exploring.
\vskip .5cm

The rest of the present section will be devoted to the discussion of  special functions determined by integral equations.
Returning  to Example \ref{exa:freud-weight}, we
 examine  the model problem for $L$. To facilitate  our discussion,  instead of the behavior \eqref{L-behavior-origin} of $L(s)$ at the origin, we      assume
\begin{equation}\label{L-behavior-origin-new}
L(s)=O\left ( \epsilon_\beta(s)  \right )~~\mbox{as}~~s\to 0,~~\epsilon_\beta(s):= \left\{ \begin{array}{ll}
1, &   1 /2 <\beta<1,\\
\ln |s|, & \beta= 1/ 2, \\
|s|^{\beta-  1 /2}, & 0<\beta< 1/2.          \end{array}
\right.
\end{equation}

The objective   is to introduce a pair of   special functions, $u_\beta$ and $v_\beta$,    as solutions of two  scalar linear integral equations.
Based on  these functions, the solution to the
 model problem $L$ of Kriecherbauer and McLaughlin \cite{Kriecherbauer-McLaughlin-1999}     can readily be constructed.
More precisely, we
assume that
$u_\beta$ and $v_\beta$ solve, respectively, the following
 integral equations
\begin{align}
&u(x)=1+\int^\infty_0 \frac {K(t) u(t) dt}{t+x},~ ~x\in (0, \infty),\label{u-integral-equation} \\
&v(x)=1-\int^\infty_0 \frac{ K(t) v(t) dt}{t+x},~ ~x\in (0, \infty),\label{v-integral-equation}
\end{align}
where
\begin{align}
          K(t) &=   \frac 1 {2\pi i} \left [
          \exp \left ( e^{(\frac \pi 2 +\frac
{\beta \pi} 2) i} t^\beta\right )
         - \exp \left ( e^{-(\frac \pi 2 +\frac {\beta
\pi} 2) i} t^\beta\right )
\right ] \label{kernel-for-continuation}  \\[.3cm]
           & =    \frac {1} \pi \exp \left (- t^\beta \sin\frac {\pi\beta} 2\right
)\sin \left ( t^\beta\cos\frac{\pi\beta} 2  \right ).\label{kernel-real}
        \end{align}
Note that
  $\eta_L(t)=2\pi i (-1)^{n+1} K(t)$ for $t\in [0, +\infty)$.

First, we formally  derive the integral equations from the model problem for $L$ satisfying \eqref{L-jump}, \eqref{L-behavior-infty} and \eqref{L-behavior-origin-new}.
From the jump condition \eqref{L-jump} and the behavior in \eqref{L-behavior-infty} and \eqref{L-behavior-origin-new}, it is readily seen  that the $(1,1)$-entry $L_{11}(s)$ is analytic in   $ \mathbb{C}\setminus (-\infty, 0]$, such that
\begin{align*}
&\left (L_{11}\right )_+(s)-\left (L_{11}\right )_-(s)=\eta_L(-s)L_{12}(s), ~s\in (-\infty, 0), \\
&L_{11}(s)=1+o(1), ~s\to\infty,  \\
&L_{11}(s)=O(\epsilon_\beta(s)), ~s\to 0.
\end{align*}
Hence, in view of the behavior of $L_{12}$ obtained from   \eqref{L-behavior-infty} and \eqref{L-behavior-origin-new},  we have from the
Plemelj formula
\begin{equation}\label{L11-integral}
L_{11}(s)=1+\frac 1 {2\pi i}  \int^0_{-\infty} \frac { \eta_L(-\tau)L_{12}(\tau) d\tau}{\tau-s}= 1+\frac 1 {2\pi i}  \int_0^{\infty} \frac { \eta_L(\tau)L_{12}(-\tau) d\tau}{-\tau-s}
\end{equation}for $s\in    \mathbb{C}\setminus (-\infty, 0]$, especially for $s\in (0, \infty)$.
Similarly,   it is also seen from  \eqref{L-jump} that $L_{12}(s)$ is analytic in   $ \mathbb{C}\setminus [0, \infty)$, and solves the scalar Riemann-Hilbert problem
\begin{align*}
&\left (L_{12}\right )_+(s)-\left (L_{12}\right )_-(s)=-\eta_L(s)L_{11}(s), ~s\in (0, \infty), \\
&L_{12}(s)=o(1), ~s\to\infty,  \\
&L_{12}(s)=O(\epsilon_\beta(s)), ~s\to 0.
\end{align*}
Hence,  we have
\begin{equation} \label{L12-integral}
L_{12}(s)=\frac 1 {2\pi i}  \int_0^{\infty} \frac { -\eta_L(\tau)L_{11}(\tau) d\tau}{\tau-s}~~\mbox{for}~~s\in    \mathbb{C}\setminus [0, \infty),
\end{equation}especially for $s\in (-\infty, 0)$. For $x\in (0, +\infty)$,
now we define
\begin{equation}\label{u-beta-v-beta}
u_\beta(x)=L_{11}(x)+(-1)^n L_{12}(-x)  ,~~v_\beta(x)=L_{11}(x)-(-1)^n L_{12}(-x),
\end{equation}$n$ being the integer appearing in the definition of $\eta_L(s)$. From \eqref{L11-integral}-\eqref{L12-integral} and \eqref{u-beta-v-beta}, it is readily verified that
$u_\beta$ and $v_\beta$ solve the
 integral equations \eqref{u-integral-equation} and \eqref{v-integral-equation}, respectively.

  It is worth noting that similar coupled scalar integral equations have   been derived in \cite{Deano-Kuijlaars-Roman-2016}, although the equations in that paper were used only to construct a contraction mapping.
Also, it is mentioned in
Fokas et al.\;\cite[p.\;161]{Fokas-Its-Kapaev-Novokshenov-2006}
that a function $u(x)$ can be parametrized via the solution of a linear integral equation. This is exactly what we are doing: We are rigorously defining a pair of new special functions,  using the above integral equations, to construct a parametrix.

We now outline the investigation carried out in Wong and Zhao \cite{Wong-Zhao-2016}. To begin with, we have
\begin{lem}\label{lem:compact-operator}
   The operator $T$ is a compact operator on $L^2[0, +\infty)$, where
   \begin{equation}\label{compact-T}
    (Tu)(x)=\int^\infty_0\frac {K(t)u(t) dt}{t+x}.
   \end{equation}
\end{lem}\vskip .2cm

  The reason is that $\frac {K(t)}{t+x}\in L^2{([0,+\infty)\times [0,+\infty))}$, with $K(t)$ given by \eqref{kernel-real}. Thus, $T$ is a Hilbert-Schmidt operator, and hence  is compact; cf. e.g. \cite[p.\;277]{Yosida-1995}. \vskip.4cm

The next step is to show that the operators $I\pm T$ has trivial null space
by using a vanishing lemma technique. To this aim, we need to know more about the analytic structure of the solutions to \eqref{u-integral-equation} and \eqref{v-integral-equation}.
It is readily seen that $(Tu)(z)$ in \eqref{compact-T} is a Stieltjes transform, and is hence analytic in the cut plane $|\arg z|<\pi$.
For later use, we need the  following  behavior of $(Tu)(z)$ at the origin:
\begin{lem}\label{lem:behavior-at-origin}
   For $u\in L^2[0, +\infty)$, we have
    $(Tu)(z)=O(\epsilon_\beta(z))$ for $|\arg z |\leq \pi$  as $z\to 0$,
  where
   $\epsilon_\beta(z)$ is given by \eqref{L-behavior-origin-new}.
\end{lem}\vskip .2cm
A description of the null space is as follows.
\begin{lem}\label{lem:vanishing-lemma}
 Assume that   $u\in  L^2[0,+\infty)$ and  $u-Tu=0$ (or $u+Tu=0$) in  $L^2[0,+\infty)$. Then,  $u\equiv 0$ for $x\in [0, \infty)$ in  $L^2[0,+\infty)$.
\end{lem}\vskip .2cm

To conclude, we have the unique existence of $u_\beta(x)$ and $v_\beta(x)$ in $L^2$ sense.
 \begin{thm}\label{thm:unique-solvability}(Wong-Zhao \cite{Wong-Zhao-2016})  There exist  unique solutions $u(x)$ and $v(x)$  to the integral equations \eqref{u-integral-equation} and \eqref{v-integral-equation}, respectively,  such that $u(x)-1\in L^2[0, +\infty)$ and $v(x)-1\in L^2[0, +\infty)$.
\end{thm}

Further analysis of the integral equations \eqref{u-integral-equation} and \eqref{v-integral-equation}  show that {$u_\beta(x)$ and $v_\beta(x)$ are in fact  bounded
 for all $\beta\in (0, 1)$.

  \begin{thm}\label{thm:boundedness} (Wong-Zhao \cite{Wong-Zhao-2016}) Let $u_\beta(x)$ and $v_\beta(x)$ be  solutions of the integral equations in \eqref{u-integral-equation} and \eqref{v-integral-equation}, respectively.
  Then
\begin{equation*}
 u_\beta(x),\;v_\beta(x)\in L^\infty[0, \infty)~~~~\mbox{and}~~~~u_\beta(x), \;v_\beta(x)\in C[0, \infty) .
\end{equation*}
\end{thm}

Since $u_\beta(x)$ and $v_\beta(x)$  solve the integral equations
\eqref{u-integral-equation} and \eqref{v-integral-equation}, with $u_\beta-1,\; v_\beta-1\in L^2[0, \infty)$, we can deduce
 that $u_\beta(z),\; v_\beta(z)=1+O(1/z)$ as $z\to\infty$ for $|\arg z|<3\pi/2$ and are of the order  $O(\epsilon_\beta(z))$ as $z\to 0$ for $|\arg z|\leq\pi$. 
Thus   we have

 \begin{thm}\label{thm: construction-of-L}(Wong-Zhao \cite{Wong-Zhao-2016}) The piece-wise analytic function
 \begin{equation}\label{model-problem-solution} L(s)=
\left\{
\begin{array}{lr}
\left(
    \begin{array}{cc}
    L_\beta(s )        &  (-1)^n U_\beta(se^{ - {\pi}   i}) \\
       (-1)^{n} U_\beta(s )  &  L_\beta(se^{-  {\pi}  i}) \\
    \end{array}
  \right), &   0<\arg s < \pi  ;
\\
  &\\
  \left(
    \begin{array}{cc}
     L_\beta(s)    & (-1)^n U_\beta(se^{ \pi   i})\\
(-1)^n U_\beta(s)  & L_\beta(se^{  \pi   i})  \\
    \end{array}
  \right), & -\pi< \arg s < 0  \end{array}
  \right .
\end{equation}
solves the Riemann-Hilbert problem \eqref{L-jump}, \eqref{L-behavior-infty} and \eqref{L-behavior-origin-new} for $L$,
where
\begin{equation}\label{L-beta-U-beta}
L_\beta(z)=\frac 1 2 \left ( u_\beta(z)+v_\beta(z)\right ),~~U_\beta(z)=\frac 1 2 \left ( u_\beta(z)-v_\beta(z)\right )
 \end{equation}for $\arg z \in (-\infty, \infty)$, $z\not=0$.
 \end{thm}

Regarding $u_\beta(z)$ and $v_\beta(z)$ as a pair of special functions,
a very
 preliminary analysis has been carried out in \cite{Wong-Zhao-2016}. For example,  for  analytic continuation, applying the Plemelj formula to the integral equations \eqref{u-integral-equation}-\eqref{v-integral-equation}, we have
\begin{equation}\label{connection-u-v} u(ze^{-\pi i})-u(ze^{\pi i})= 2\pi i K(z) u(z),~~v(ze^{-\pi i})-v(ze^{\pi i})= -2\pi i K(z) v(z),\end{equation}
initially for real and positive $z$, and then for complex  $z$ with $\arg z\in (-\infty, \infty)$; cf. \eqref{kernel-for-continuation}.

Also, treating the  integral equations as Stieltjes transform   and applying Theorem \ref{thm:unique-solvability},
we readily obtain the asymptotic approximations
\begin{equation}\label{u-beta-v-beta-infty}
 u_\beta (z)\sim 1+\sum^\infty_{k=1} \frac {c_k}{z^k} ~~\mbox{and}~~ v_\beta(z)\sim 1+\sum^\infty_{k=1} \frac {d_k}{z^k}~~~~\mbox{as}~ z \to+\infty,~~|\arg z| <   \frac {3\pi} 2,
\end{equation}
where
\begin{equation}\label{u-beta-v-beta-infty-coeff}
  c_k=(-1)^{k-1} \int ^\infty_0 K(t) t^{k-1}u_\beta(t)dt~~\mbox{and}~~d_k=(-1)^{k} \int ^\infty_0 K(t) t^{k-1}v_\beta(t)dt
\end{equation}for $k=1,2,\cdots$.

The Stokes phenomenon is also addressed in \cite{Wong-Zhao-2016}.  Assume that the functions $u_\beta(z)$ and $v_\beta(z)$  solve respectively \eqref{u-integral-equation} and \eqref{v-integral-equation}, such that $u_\beta-1,~v_\beta-1\in   L^2[0, \infty)$. Then the Stokes lines for both functions are
\begin{equation*}
  \arg z=\pm \left\{  \frac \pi 2(3-\alpha)+2(l-1)\pi\right\},
\end{equation*}
  where  $\alpha=\frac 1 \beta \in ( 4l-3, 4l+1]$, $l=1,2,\cdots$.

Much remains to be done.  We propose a thorough  investigation of the analytic and asymptotic  properties of the  functions  $u_\beta$ and $v_\beta$,   such as
  their  zeros,  modulus function  \cite[\S 10.18]{nist-2010}, kernel  \cite{Kuijlaars-Vanlessen-2002}, and differential or difference  equations, etc. Here, we list a few.
    \begin{description}
    \item (i) Determine   the coefficients $c_k$ and $d_k$ in the expansions at infinity given by \eqref{u-beta-v-beta-infty}-\eqref{u-beta-v-beta-infty-coeff}, which  are now expressed in terms of $u_\beta$ and $v_\beta$. It would be interesting and challenging to decode the coefficients
from \eqref{u-integral-equation} and \eqref{v-integral-equation} in an explicit manner.
    \item (ii)  As a refinement  of Theorem \ref{thm:boundedness}, one might  reasonably expect that the behaviors of $u_\beta$ and $v_\beta$ at the origin are of the form $\sum_{k=0}\tilde c_k z^k+\sum_{k=1} \tilde d_k z^{\beta k}$. A natural problem is to determine the small-$z$ asymptotic formula and to evaluate   explicitly the coefficients  $\tilde c_k$ and  $\tilde d_k$, at least the leading coefficients such as  $k=0$ and $k=1$.
\item (iii) The eigenvalues of the compact operator $T$ given by \eqref{compact-T}, that is, find $\lambda$ such that
the null space of  $\lambda I-T$ is nontrivial in $L^2[0, +\infty)$.
\item (iv) It has been conjectured in \cite{Wong-Zhao-2016} that the function $u_\beta(x)$ is monotonically decreasing in $x\in [0, \infty)$, and also that the initial value $u_\beta(0)$ is decreasing in $\beta\in (0, 1)$.
Corresponding conjectures can be made for the function $v_\beta(x)$, such as
it is increasing in $x\in [0, \infty)$, and also that the initial value $v_\beta(0)$ is increasing in $\beta\in (0, 1)$.
It is also of  interest to consider the complete monotonicity of these special functions.
  \end{description}


\section{Wiener-Hopf equation}\label{sec: integral-eq}
As far as we are aware, the first paper that addresses asymptotic solutions to an integral equation of the form in \eqref{int-eq} is by Muki and Sternberg \cite{Muki-Sternberg-1970} in 1970. But, their approach assumes that the solution behaves like $t^{-\delta}/M$ as $t\to+\infty$, where $\delta$ is a positive number and $M$ is a constant.  An attempt was made by Li and Wong \cite{Li-Wong-1994} in 1994 to remove such an assumption and to derive the asymptotic solution directly from the integral equation, but their attempt was not successful. Later, it was realized that equation \eqref{int-eq} is known as of Wiener-Hopf type. This naturally led Li and Wong \cite{Li-Wong-2021} to use the famous Wiener-Hopf technique \cite{Wiener-Hopf-1931}.  Here we first give an outline of this technique, leading to an integral representation of the solution. For convenience we reproduce \eqref{int-eq} below;
\begin{equation}\label{int-eq-WH}
u(t)=f(t)+\int^\infty_0 k(t-\tau) u(\tau) d\tau, ~~~~t>0.
\end{equation}

 The following notations will be used throughout this section:
 \begin{equation}\label{u-plus}
 u_+(t)=\left\{\begin{array}{ll}
                u(t), & t>0, \\[.3cm]
                 0, & t<0
               \end{array}
        \right .
 \end{equation}and
 \begin{equation}\label{u-minus}
 u_-(t)=\left\{\begin{array}{ll}
               0, & t>0, \\[.3cm]
             f(t)+\int^\infty_0k(t-\tau) u(\tau)d\tau,  & t<0.
               \end{array}
        \right .
 \end{equation}
 Since we are only concerned with $t>0$ in equation \eqref{int-eq-WH}, in most cases the function $f(t)$  is identically zero for $t$ negative. Equation
 \eqref{int-eq-WH} can now  be written as
 \begin{equation}\label{int-eq-WH-modify}
u_+(t)+u_-(t)=f(t)+\int^\infty_{-\infty} k(t-\tau) u_+(\tau) d\tau, ~~~~-\infty<t<\infty.
\end{equation}
 Taking Fourier transforms on both sides of this equation gives
\begin{equation}\label{WH-eq-Fourier}
 \left [ 1-\sqrt{2\pi}K(\lambda)\right ] U_+(\lambda)=F(\lambda)-U_-(\lambda),
\end{equation}
 where
 \begin{equation}\label{Fourier-trans}
 \begin{array}{l}
  \displaystyle{K(\lambda)=\frac 1 {\sqrt{2\pi}} \int^\infty_{-\infty} k(t) e^{i\lambda t}dt,~~~~~~
 F(\lambda)=\frac 1 {\sqrt{2\pi}} \int^\infty_{-\infty} f(t) e^{i\lambda t}dt,} \\[.4cm]
  \displaystyle{U_\pm(\lambda)=\frac 1 {\sqrt{2\pi}} \int^\infty_{-\infty} u_\pm(t) e^{i\lambda t}dt,}
 \end{array}
 \end{equation}
    and $\lambda=\sigma+i\tau$.  In \eqref{WH-eq-Fourier}, we note that
$U_+(\lambda)$ and $U_-(\lambda)$ are unknown functions, and that $K(\lambda)$ and $F(\lambda)$ are known functions. To solve this functional equation involves a factorization of the function $1-\sqrt{2\pi}K(\lambda)$, which is the crux of the idea in the Wiener-Hopf technique.

Let $L(\lambda)$ be an analytic function of  $\lambda=\sigma+i\tau$ in the strip $\tau_-<\tau<\tau_+$ such that $\left |L(\sigma+i\tau)\right |\leq C|\sigma|^{-p}$, $p>0$ and $C>0$, as $|\sigma|\to\infty$. This inequality holds uniformly for $\tau$ in the strip $\tau_-+\varepsilon\leq\tau\leq\tau_+-\varepsilon$, $\varepsilon>0$. For
 $\tau_-<c<\tau<d<\tau_+$, we have by Cauchy's theorem
   \begin{equation*}
    L(\lambda)=\frac 1 {2\pi i} \int_\mathcal{C}\frac{L(\zeta)}{\zeta-\lambda} d\zeta,
   \end{equation*}
where     $\mathcal{C}$, a positively oriented contour, is the boundary of the rectangle whose four corners are at $-R+ic$,  $R+ic$, $R+id$ and  $-R+id$.
Letting $R\to +\infty$ yields the decomposition
\begin{equation}\label{L-decomposition}
L(\lambda)=L_+(\lambda)-L_-(\lambda) ,
\end{equation}where
  \begin{equation}\label{L-bv}
   L_+(\lambda)=\frac 1 {2\pi i} \int_{-\infty+ic}^{\infty+ic}\frac{L(\zeta)}{\zeta-\lambda} d\zeta,
  ~~~~
    L_-(\lambda)=\frac 1 {2\pi i} \int_{-\infty+id}^{\infty+id}\frac{L(\zeta)}{\zeta-\lambda} d\zeta,
  \end{equation}
    and $\tau_-<c <d<\tau_+$. It is easily seen that
 $L_+(\lambda)$ is analytic for all $\tau>c$ and $L_-(\lambda)$ is analytic for all $\tau<d$.

 If $\ln K(\lambda)$ satisfies the conditions imposed on the function $L(\lambda)$ above, then we immediate have
 \begin{equation}\label{K-quotient}
 K(\lambda)=\frac {K_+(\lambda)}{K_-(\lambda)},
 \end{equation}where
  \begin{equation}\label{K-plus}
   K_+(\lambda)=\exp\left\{\frac 1 {2\pi i} \int_{-\infty+ic}^{\infty+ic}\frac{\ln K(\zeta)}{\zeta-\lambda} d\zeta\right\}
  \end{equation}and
 \begin{equation}\label{K-minus}
   K_-(\lambda)=\exp\left\{\frac 1 {2\pi i} \int_{-\infty+id}^{\infty+id}\frac{\ln K(\zeta)}{\zeta-\lambda} d\zeta\right\}.
  \end{equation}
It is not difficult to see that  $K_+(\lambda)$ and $K_-(\lambda)$ are analytic, bounded and nonzero in the half-planes $\tau>\tau_-$ and $\tau<\tau_+$, respectively.

Returning to \eqref{WH-eq-Fourier}, we  assume that $1-\sqrt{2\pi}K(\lambda)$ has zeros $a_1,~\cdots~,~a_m$. As in \eqref{K-quotient}, we can find
$K_+(\lambda)$ and $K_-(\lambda)$, which are analytic  and  zero-free in the respective half-plane  $\tau>\tau_-$ and $\tau<\tau_+$, such that
\begin{equation}\label{K-comb}
1-\sqrt{2\pi}K(\lambda)
=\frac {K_+(\lambda)}{K_-(\lambda)}(\lambda-a_1)(\lambda-a_2)\cdots(\lambda-a_m).
\end{equation}
Equation \eqref{WH-eq-Fourier} may now be written as
\begin{equation}\label{WH-eq-Fourier-new}
U_+(\lambda) K_+(\lambda) (\lambda-a_1)(\lambda-a_2)\cdots(\lambda-a_m)
=K_-(\lambda)F(\lambda)-K_-(\lambda)U_-(\lambda).
\end{equation}
As in \eqref{L-decomposition}, the term $K_-(\lambda)F(\lambda)$ can be decomposed into the form
\begin{equation}\label{KF-decomposite}
K_-(\lambda)F(\lambda)=C_+(\lambda)+C_-(\lambda),
\end{equation}where $C_+(\lambda)$ and $C_-(\lambda)$ are analytic in
 $\tau>\tau_-$ and $\tau<\tau_+$, respectively. With \eqref{KF-decomposite},
 we rearrange \eqref{WH-eq-Fourier-new} so that we can define a function $E(\lambda)$ by
\begin{equation}\label{E-def}
E(\lambda)= U_+(\lambda) K_+(\lambda) (\lambda-a_1)(\lambda-a_2)\cdots(\lambda-a_m)-C_+(\lambda)
=C_-(\lambda) -K_-(\lambda)U_-(\lambda).
\end{equation}
 This equation defines $E(\lambda)$ only in the strip $\tau_-<\tau<\tau_+$. But, by the first equality $E(\lambda)$ is actually defined and analytic in
 $\tau>\tau_-$, and by the second equality $E(\lambda)$ is   defined and analytic in $\tau<\tau_+$. Hence, using analytic continuation, we can define $E(\lambda)$ in the whole $\lambda$-plane. Suppose that it can be shown
\begin{equation*}
 \left | U_+(\lambda) K_+(\lambda) (\lambda-a_1)(\lambda-a_2)\cdots(\lambda-a_m)-C_+(\lambda)\right |\leq |\lambda|^p, ~~~~\mbox{as}~\lambda\to\infty~\mbox{in}~  \tau>\tau_-,
\end{equation*}
\begin{equation*}
\left |C_-(\lambda) -K_-(\lambda)U_-(\lambda)\right |\leq |\lambda |^q,~~~~\mbox{as}~\lambda\to\infty~\mbox{in}~  \tau<\tau_+.
\end{equation*}
Then, by the extended form of the Liouville theorem, $E(\lambda)$ is a polynomial $P(\lambda)$ of degree less than or equal to the integer part of $\min\{p, q\}$, i.e.,
\begin{equation}\label{U-in-polynomial}
\begin{array}{r}
 U_+(\lambda) K_+(\lambda) (\lambda-a_1)(\lambda-a_2)\cdots(\lambda-a_m)-C_+(\lambda)=P(\lambda),\\[.4cm]
 C_-(\lambda) -K_-(\lambda)U_-(\lambda)=P(\lambda).
\end{array}
\end{equation}
 To determine the degrees $p$ and $q$, we need to know the behaviors of $U_+(\lambda)$ and $U_-(\lambda)$. In special cases, these can be found from the behavior of $u(t)$. For example, if $u(t)\sim t^{-1/2}$ as $t\to 0^+$, then by the Abelian theorem $U_+(\lambda)\sim \lambda^{-\frac 1 2}$ as $\lambda\to\infty$ in the upper half-plane $\tau>0$; see \cite[p.\;36 and p.\;69]{Noble-1958}.  In equations in \eqref{U-in-polynomial},
 $U_+(\lambda)$ and $U_-(\lambda)$ are determined only within an arbitrary polynomial $P(\lambda)$; that is, within a finite number of arbitrary constants, which must be determined by some other methods.

 By
using the analytic properties of the function $U_+(\lambda)$ and Cauchy's theorem, we first derive the integral representation
\begin{equation}\label{U-plus-representation}
U_+(\lambda)=\frac 1 {2\pi i} \int_{-\infty+i\alpha}^{\infty+i\alpha}\frac{U_+(\zeta)}{\zeta-\lambda} d\zeta,
\end{equation}
where $\alpha$ is a real number bigger than $\tau_-$ and $\lambda$ is any point in the upper half-plane $\tau>\alpha$. The denominator of the above Cauchy integral can be written as
\begin{equation}\label{denominator}
\frac 1 {i(\zeta-\lambda)}=\int^\infty_0 e^{-i (\zeta-\lambda)t}dt,
\end{equation}
where $\lambda$ is any point in the half-plane  $\tau>\alpha$.
Substituting \eqref{denominator} into
\eqref{U-plus-representation} and interchanging the order of integration
give
\begin{equation*}
U_+(\lambda)=\frac 1 {\sqrt{2\pi}}
 \int_0^{\infty}
 e^{i \lambda t}\left [\frac 1 {\sqrt{2\pi}} \int_{-\infty+i\alpha}^{\infty+i\alpha} U_+(\zeta)e^{-i t\zeta }  d\zeta
\right ] dt.
\end{equation*}
Comparing the integral on the right side of this equation with the Fourier transforms given in \eqref{Fourier-trans}, we have
\begin{equation}\label{u-inverse-Fourier}
u(t)=\frac 1 {\sqrt{2\pi}} \int_{-\infty+i\alpha}^{\infty+i\alpha} U_+(\zeta)e^{-i t\zeta }  d\zeta,~~~~t>0.
\end{equation}
On account of \eqref{U-in-polynomial}, the solution of \eqref{int-eq-WH} has the integral representation
\begin{equation}\label{u-in-P}
u(t)=\frac 1 {\sqrt{2\pi}} \int_{-\infty+i\alpha}^{\infty+i\alpha}
\frac {P(\lambda)+C_+(\lambda)}
{K_+(\lambda) (\lambda-a_1)(\lambda-a_2)\cdots(\lambda-a_m)}\;
 e^{-i t\lambda }  d\lambda,~~~~t>0,
\end{equation}for any $\alpha$, $\tau_-<\alpha<\tau_+$.

A major assumption in the preceding analysis is that the kernel function $k(t)$ in equation \eqref{int-eq-WH} is exponentially decaying at infinity.
This enables one to factor
 \begin{equation*}
 \frac { 1-\sqrt{2\pi}K(\lambda)}{ (\lambda-a_1)(\lambda-a_2)\cdots(\lambda-a_m)}
 \end{equation*}
 into two components
 $K_+(\lambda)$ and $K_-(\lambda)$, one analytic in an upper half-plane, the other analytic in a lower half-plane, and the two half-planes overlapping in  an infinite strip parallel to the real axis.

 If the kernel is not exponentially decaying, e.g., the Cauchy density function
\begin{equation*}
 k(t)=\frac 1 \pi \frac 1 {1+t^2},
\end{equation*}
then one would not to be able to decompose
$\left[ 1-\sqrt{2\pi}K(\lambda)\right ]/ (\lambda-a_1)(\lambda-a_2)\cdots(\lambda-a_m)$
into two factors analytic in two overlapping half-planes. A  natural
way to extend the method presented above is to introduce an exponential  function, such as $e^{-\varepsilon |t|}$, in the kernel.  This idea has been used by Carlson and Heins \cite{Carlson-Heins-1947} and by Carrier \cite{Carrier-1956}, but the  arguments in both papers are only formal. Rigorous justifications were given by Widom \cite{Widom-1958} for the case $k \in L^2 (\mathbb{R})$,  and by Krein for $k \in L^1 (\mathbb{R})$. (For a brief introduction to the work of Krein \cite{Krein-1958} and  Gohberg and Krein \cite{Gohberg-Krein-1958}, see \cite[Ch.\;20]{Beals-Wong-2020}.)
Here, we give a quick review of the argument given in Widom \cite{Widom-1958}.

Let $\{\varepsilon_n\}$ be a monotonically decreasing sequence tending to zero as $n\to\infty$. Define
\begin{equation*}
\tilde e_n(t)=\left\{
\begin{array}{ll}
 e^{-\varepsilon_n t}, & t>0 \\[.3cm]
  0, & t<0.
\end{array}\right .
\end{equation*}
We call $u(t)$ a solution of \eqref{int-eq-WH}, if for any $\eta>0$ it satisfies
\begin{equation}\label{weak-sol}
\lim_{n\to \infty}\left \|
e^{-\eta |t|}\int_{-\infty}^\infty k(t-\tau) \tilde e_n(\tau) u(\tau)d\tau
-e^{-\eta |t|}\left [u(t)-f(t)\right ]\right \|_2=0.
\end{equation}
In \eqref{weak-sol}, $e^{-\eta |t|}$ is inserted for some technical reasons. If such a solution does exist, then it is well known that there will exist a subsequence
$\{\varepsilon_{n_k}\}$ such that $\varepsilon_{n_k}\to 0$ as $k\to\infty$ and \begin{equation}\label{weak-int-eq}
\lim_{k\to \infty}
 \int_0^\infty k(t-\tau) \tilde e_{n_k}(\tau) u(\tau)d\tau
=u(t)-f(t)
\end{equation}
holds almost everywhere; see \cite[Prop.\;0.1.10]{Butzer-Nessel-1971}. To make this argument work, Widom first proved that the functional equation
\begin{equation}\label{WH-eq-Fourier-md}
U_-(\zeta)=F(\zeta)+\left [\sqrt{2\pi}K(\zeta) -1\right ] U_+(\zeta),~~~~\zeta\in \mathbb{R},
\end{equation}holds almost everywhere; cf.
\eqref{WH-eq-Fourier}. Next, he introduced the function
\begin{equation}\label{Psi-def}
\Psi(\zeta):=\frac {\left[1-\sqrt{2\pi} K(\zeta)\right ](\zeta+i)^{\frac  n  2}
(\zeta-i)^{\frac n  2}}{(\zeta-\alpha_1)^{m_1}(\zeta-\alpha_2)^{m_2}\cdots
(\zeta-\alpha_l)^{m_l}}\left (\frac {\zeta-i}{\zeta+i}\right )^k,
\end{equation}
which has no zero on the real line. Here, $\alpha_j$ for $1\leq j\leq l$ is a zero of
$1-\sqrt{2\pi} K(\zeta)$ with multiplicity $m_j$, and
$n=m_1+\cdots+m_l$. To obtain a factorization formula like the one in \eqref{K-comb}, he set
\begin{equation}\label{Cauchy-integral}
 \chi(\lambda):= \frac 1 { 2\pi i} \int_{-\infty }^{\infty }
\frac {\ln \Psi(\zeta)}{\zeta-\lambda}  d\zeta;
\end{equation}
cf. \eqref{K-plus} and \eqref{K-minus}. Note that this integral defines an analytic function in $\mathbb{C}\setminus \mathbb{R}$. Put
\begin{equation}\label{phi-in-chi}
\begin{array}{ll}
   \phi_+(\lambda):=(\lambda+i)^{\frac  n 2-k} e^{-\chi_+(\lambda)} , & \Im \lambda >0, \\[.3cm]
 \phi_-(\lambda):=(\lambda-i)^{-\frac  n 2-k} e^{-\chi_-(\lambda)}, & \Im \lambda <0,
\end{array}
\end{equation}
where
\begin{equation*}
\chi_+(\lambda)=\chi(\lambda)~~~~\mbox{in}~~\Im \lambda>0,~~~~
\chi_-(\lambda)=\chi(\lambda)~~~~\mbox{in}~~\Im \lambda<0
\end{equation*}with boundary values
\begin{equation}\label{Cauchy-BV}
 \begin{array}{l}
 \displaystyle{\chi_+(\zeta)=\lim_{\varepsilon\to 0^+} \chi(\zeta+i\varepsilon )
 = \frac 1 2 \ln \Psi(\zeta)+   \frac 1 { 2\pi i}
 \int_{-\infty}^{\infty}{\!\!\!\!\!\!\!\!\!\!\!\!\!    {}-{} \,\,\,\,\,\, }
 \frac {\ln \Psi(\xi)}{\xi-\zeta }  d\xi,} \\[.4cm]
\displaystyle{\chi_-(\zeta)=\lim_{\varepsilon\to 0^+} \chi(\zeta-i\varepsilon )
 = -\frac 1 2 \ln \Psi(\zeta)+   \frac 1 { 2\pi i} \int_{-\infty}^{\infty}{\!\!\!\!\!\!\!\!\!\!\!\!\!    {}-{} \,\,\,\,\,\, }
\frac {\ln \Psi(\xi)}{\xi-\zeta }  d\xi }
 \end{array}
\end{equation}for $\zeta\in \mathbb{R}$, where the  bar indicates that the integral is a Cauchy principal value.
It is now readily seen that
\begin{equation}\label{K-comb-new}
1-\sqrt{2\pi}K(\zeta)=\frac {\phi_-(\zeta)}{\phi_+(\zeta)}(\zeta-\alpha_1)^{m_1}(\zeta-\alpha_2)^{m_2}
\cdots (\zeta-\alpha_l)^{m_l};
\end{equation}
cf. \eqref{K-comb}. Applying \eqref{K-comb-new} to the functional equation \eqref{WH-eq-Fourier-md} yields
\begin{equation}\label{U-plus-minus}
\frac{U_+(\zeta)}{\phi_+(\zeta)}(\zeta-\alpha_1)^{m_1}(\zeta-\alpha_2)^{m_2}
\cdots (\zeta-\alpha_l)^{m_l}
=-\frac{U_-(\zeta)} {\phi_-(\zeta)}
+\frac{F(\zeta)}{\phi_-(\zeta)}.
\end{equation}
To obtain an analogue of \eqref{KF-decomposite}, Widom set
\begin{equation}\label{Cauchy-integral-H}
 H(\lambda) = \frac 1 { 2\pi i} \int_{-\infty }^{\infty }
\frac {F(\zeta)}{\phi_-(\zeta)}\frac {d\zeta}{\zeta-\lambda}
\end{equation}with
\begin{equation*}
H_+(\lambda)=H(\lambda)~~~~\mbox{in}~~\Im \lambda>0,~~~~
H_-(\lambda)=H(\lambda)~~~~\mbox{in}~~\Im \lambda<0
\end{equation*}and
\begin{equation}\label{Cauchy-BV-H}
 \begin{array}{l}
 \displaystyle{H_+(\zeta)=
 \frac 1 2 \frac {F(\zeta)}{\phi_-(\zeta)}+   \frac 1 { 2\pi i}
 \int_{-\infty}^{\infty}{\!\!\!\!\!\!\!\!\!\!\!\!\!    {}-{} \,\,\,\,\,\, }
\frac {F(\xi)}{\phi_-(\xi)}\frac {d\xi}{\xi-\zeta},} \\[.4cm]
\displaystyle{H_-(\zeta)=
 -\frac 1 2 \frac {F(\zeta)}{\phi_-(\zeta)}+   \frac 1 { 2\pi i}
\int_{-\infty}^{\infty}{\!\!\!\!\!\!\!\!\!\!\!\!\!    {}-{} \,\,\,\,\,\, }
\frac {F(\xi)}{\phi_-(\xi)}\frac {d\xi}{\xi-\zeta}}
 \end{array}
\end{equation}for $\zeta\in \mathbb{R}$; cf. \eqref{Cauchy-BV}. From \eqref{Cauchy-BV-H}, it follows
\begin{equation}\label{F-phi-decomposition}
\frac{F(\zeta)}{\phi_-(\zeta)}=H_+(\zeta)-H_-(\zeta).
\end{equation}
Coupling \eqref{U-plus-minus} and \eqref{F-phi-decomposition}
yields
\begin{equation}\label{U-plus-minus-net}
\frac{U_+(\zeta)}{\phi_+(\zeta)}(\zeta-\alpha_1)^{m_1}(\zeta-\alpha_2)^{m_2}
\cdots (\zeta-\alpha_l)^{m_l}-H_+(\zeta)
=-\frac{U_-(\zeta)} {\phi_-(\zeta)}
-H_-(\zeta).
\end{equation}
The left-hand side of \eqref{U-plus-minus-net}
defines an analytic function on the upper half-plane $\mathbb{C}_+$, and the right-hand side defines an analytic function on the lower half-plane
$\mathbb{C}_-$; these two functions agree on the real line $\mathbb{R}$. By using a theorem of Carleman  \cite[p.\;40]{Carleman-1944}, Widom showed that these two functions are analytic continuations of each other and together represent an entire function $P(\lambda)$. Furthermore, by using an estimate of the function $\phi_-(\lambda)$ in \eqref{phi-in-chi}, Widom concluded that  $P(\lambda)$ is in fact a polynomial of degree less than $[\frac n 2+k]$. From \eqref{U-plus-minus-net}, we have
\begin{equation}\label{U-plus-rep-Widom}
U_+(\lambda)=\frac{\left [
P(\lambda)+H_+(\lambda)\right ]\phi_+(\lambda)}
{(\lambda-\alpha_1)^{m_1}(\lambda-\alpha_2)^{m_2}
\cdots (\lambda-\alpha_l)^{m_l}}.
\end{equation}
By inversion, we obtain
\begin{equation}\label{u-rep-Widom}
u(t)=
\frac 1 {\sqrt{2\pi}}\int_{-\infty+i\gamma}^{ \infty+i\gamma}
\frac{\left [
P(\lambda)+H_+(\lambda)\right ]\phi_+(\lambda)}
{(\lambda-\alpha_1)^{m_1}(\lambda-\alpha_2)^{m_2}
\cdots (\lambda-\alpha_l)^{m_l}} e^{-it\lambda} d\lambda, ~~~~t>0,
\end{equation}
for any $\gamma>0$;  cf. \eqref{u-inverse-Fourier}.  In his final step, Widom established the limit in \eqref{weak-sol}, thus proving that $u(t)$ is a solution of \eqref{int-eq-WH}.

For illustration purpose, it suffices to consider just the special case in which all zeros of $1-\sqrt{2\pi}K(\lambda)$ are simple; i.e., $m_1=\cdots=m_l=1$, and $K(\lambda)$ and $F(\lambda)$ are both smooth, except possibly for a finite number of points in $\mathbb{R}\setminus \{\alpha_1, \cdots, \alpha_l\}$.  By Cauchy's residue theorem, the integral in \eqref{u-rep-Widom} can be expressed as
\begin{equation}\label{u-rep-simple}
\begin{array}{rcl}
 u(t)&=&\displaystyle{
\frac 1 {\sqrt{2\pi}}
\int_{-\infty}^{\infty}{\!\!\!\!\!\!\!\!\!\!\!\!\!    {}-{} \,\,\,\,\,\, }
\frac{\left [
P(\lambda)+H_+(\lambda)\right ]\phi_+(\lambda)}
{(\lambda-\alpha_1) (\lambda-\alpha_2)
\cdots (\lambda-\alpha_l) } e^{-it\lambda} d\lambda}\\[.4cm]
&&\displaystyle{- i\;\sqrt{\frac \pi 2}\;\sum^l_{j=1} G_j(\alpha_j) e^{-it\alpha_j}
, }
\end{array}
\end{equation}
where
\begin{equation}\label{G-j-def}
G_j(\lambda)=\frac{\left [
P(\lambda)+H_+(\lambda)\right ]\phi_+(\lambda) (\lambda-\alpha_j)}
{(\lambda-\alpha_1) (\lambda-\alpha_2)
\cdots (\lambda-\alpha_l) }.
\end{equation}
Note that here we have made use of the formula
\begin{equation}\label{formula}
\frac 1 {\sqrt{2\pi}}
\int_{-\infty}^{\infty}{\!\!\!\!\!\!\!\!\!\!\!\!\!    {}-{} \,\,\,\,\,\, }
\frac{e^{-it\xi} }
{\xi} d\xi=
-\sqrt{\frac \pi 2}\; i.
\end{equation}

Finally, choose $l$ cut-off functions $\eta_j(\lambda)$, $j=1,\cdots, l$, such that
$\eta_j(\lambda)$ is equal to $1$ for $\lambda$ near $\alpha_j$, $j=1,\cdots, l$, and $\mbox{supp~} \eta_j \cap \mbox{supp~} \eta_i=\emptyset$ for each $j\not=i$. Let $\eta_0=1-\sum^l_{j=1} \eta_j$, the integral in \eqref{u-rep-simple} can now be written as
\begin{equation}\label{u-rep-simple-final}
\begin{array}{rcl}
 u(t)&=&\displaystyle{
\frac 1 {\sqrt{2\pi}}\sum^l_{j=1}
e^{-it\alpha_j}
\int_{-\infty }^{ \infty }
\frac{G_j(\xi)\eta_j(\xi)-G_j(\alpha_j)}
{ \xi-\alpha_j  } e^{-it(\xi-\alpha_j)}  d\xi}\\[.4cm]
&&
\displaystyle{+
\frac 1 {\sqrt{2\pi}} \int_{-\infty }^{ \infty }
\frac{\left [
P(\xi)+H_+(\xi)\right ]\phi_+(\xi)\eta_0(\xi)}
{(\xi-\alpha_1) (\xi-\alpha_2)
\cdots (\xi-\alpha_l) } e^{-it\xi} d\xi}\\[.4cm]
&&\displaystyle{-\sqrt{ 2 \pi}\; i\sum^l_{j=1} G_j(\alpha_j) e^{-it\alpha_j}
. }
\end{array}
\end{equation}

Asymptotic behavior of the solution $u(t)$ to equation \eqref{int-eq-WH}
can be obtained from the Fourier integrals in \eqref{u-rep-simple-final};
see \cite[Thm.\;1, p.\;199 and Ex.\;11, p.\;234]{Wong-1989}.
But, the amount of work to derive just the leading term in the asymptotic expansion of $u(t)$ is still tremendous. This is because the integrands in
\eqref{u-rep-simple-final} involve functions $\phi_+(\lambda)$ and $H_+(\lambda)$, which are given in terms of Hilbert transforms; see \eqref{phi-in-chi}, \eqref{Cauchy-BV} and \eqref{Cauchy-BV-H}. Fortunately, asymptotic expansions of Hilbert transforms can be obtained via corresponding results for the Stieltjes transform
\begin{equation}\label{Stieltjes-tr}
S_\varphi(z)=\int^\infty_0\frac {\varphi(t)}{t+z} dt,~~~~|\arg z|<\pi,
\end{equation}
and the one-sided Hilbert transform
\begin{equation}\label{Hilbert-one-side}
H^+_\varphi(x)=\int^\infty_0{\!\!\!\!\!\!\!\!\!\!\!\! {}-{} \,\,\,\,\,  }\frac {\varphi(t)}{t-x} dt,~~~~x\in \mathbb{R}^+,
\end{equation}
where again the bar indicates that the integral is a Cauchy principal value.
In deriving the asymptotic expansions of these transforms, we also run into the Mellin transform of $\varphi$, which is defined by
\begin{equation}\label{Mellin-tr}
M[\varphi; z]=\int^\infty_0 t^{z-1}\varphi(t) dt, ~~~~~~~~~~\Re z>0.
\end{equation}
We assume that $\varphi\in L^1(\mathbb{R}^+)$ has an asymptotic expansion of the form
\begin{equation}\label{varphi-behavior}
\varphi(t)\sim \sum^\infty_{s=0} \sum^s_{i=0} \alpha_{s, i} t^{-s-\theta} \ln^i t,~~~t\to\infty,
\end{equation}for fixed $\theta\in (0, 1]$.  Let
\begin{equation}\label{varphi-error}
\varphi_n(t)= \varphi(t)-  \sum^{n-1}_{s=0} \sum^s_{i=0} \alpha_{s, i} t^{-s-\theta} \ln^i t.
\end{equation}
For convenience, we also introduce a new notation. Recall the Pochhammer symbol $(\theta)_n=\theta(\theta+1)\cdots (\theta+n-1)$.  Since
\begin{equation*}
(-1)^n \frac {d^n}{dt^n} \left (\frac {t^{-\theta}}{(\theta)_n}\right )=t^{-n-\theta}
\end{equation*}
for any $\theta\in \mathbb{R}$ and $n\in \mathbb{N}$ with $(\theta)_n\not= 0$, we have
\begin{equation}\label{derivatives}
(-1)^{n+i} \frac {d^n}{dt^n} \left\{
\frac {d^i}{d\theta^i}
\left (\frac {t^{-\theta}}{(\theta)_n}\right )
\right\}
=t^{-n-\theta}\ln^i t
\end{equation}for any $i\in \mathbb{N}$. Let
\begin{equation}\label{H-def}
H_{\theta, n, i}(t) =(-1)^{n+i}
\frac {d^i}{d\theta^i}
\left (\frac {t^{-\theta}}{(\theta)_n}\right )
=\sum^i_{j=0} C^j_{\theta, n, i}  t^{-\theta}\ln^j t,
\end{equation}
where
\begin{equation}\label{H-coef-def}
C^j_{\theta, n, i}=(-1)^{n+i+j} \begin{pmatrix}
                                  i \\
                                  j \\
                                \end{pmatrix}
\frac {d^{i-j}}{d\theta^{i-j}}
\left (\frac {1}{(\theta)_n}\right ).
\end{equation}
Note that  from \eqref{H-def}, we have
\begin{equation*}
H_{\theta, n, i}^{(n)}(t) = t^{-n-\theta } \ln^i t.
\end{equation*}

\begin{thm}
Suppose that $\varphi$ satisfies  \eqref{varphi-behavior} with $\theta\in (0, 1)$. Then for any  $z\in \mathbb{C}$ with $|\arg z|<\pi$ and $n\in \mathbb{N}^+$, we have
\begin{equation}\label{Stieltjes-expan}
\begin{array}{l}
\displaystyle{  S_\varphi(z)= \sum^{n-1}_{s=0} \sum^s_{i=0} \sum^i_{j=0}
(-1)^j \pi \alpha_{s, i} C^j_{\theta, s, i}
\frac {d^j}{d\theta^j}\left [ \frac {(\theta)_s}{\sin\theta\pi} z^{-\theta}\right ] z^{-s}
}
\\[.5cm]
\displaystyle{
~~~~~~~~~~     -
\sum^n_{s=1} (s-1)! c_s z^{-s} +\frac {(-1)^n}{z^n} \int^\infty_0\frac {\tau^n \varphi_n(\tau)}{\tau+z} d\tau,}
\end{array}
\end{equation}
where $c_s=\frac{(-1)^s}{(s-1)!} M[\varphi; s]$ and $M[\varphi; s]$ is the Mellin transform of $\varphi$.
\end{thm}

To state the next result, we need another notation, namely
\begin{equation}\label{gamma-k-s}
\gamma_{k,s}=\int^\infty_0\frac {\ln^k t}{(t+1)^{s+2}} dt,
\end{equation}
where $s,~k\in \mathbb{N}$. Direct calculation gives
\begin{equation}\label{gamma-0-s}
\gamma_{0,s}=\frac 1 {s+1},~~~~\gamma_{1,s}=-\frac 1 {s+1} \sum^ s_{i=1} \frac 1 i ,
\end{equation}
\begin{equation}\label{gamma-2-s}
\gamma_{2,s}= \frac 4 {s+1} \sum^\infty_{i=1} \frac {(-1)^{i-1}} { i^2}  -\frac 2 {s+1} \sum^ s_{i=1}
\int^\infty_0 \frac {\ln t}{(t+1)^{i+1}} dt
\end{equation}
and
\begin{equation}\label{gamma-k-s-formula}
\begin{array}{rl}
 \gamma_{k,s} & \displaystyle{= \frac {[1+(-1)^k] k!} {s+1} \sum^\infty_{i=1} \frac {(-1)^{i-1}} { i^k}  -\frac k {s+1} \sum^s_{i=1}
\int^\infty_0 \frac {\ln^{k-1} t}{(t+1)^{i+1}} dt} \\[.5cm]
    &   \displaystyle{=\frac {[1+(-1)^k] k!} {s+1} \sum^\infty_{i=1} \frac {(-1)^{i-1}} { i^k}  -\frac k {s+1} \sum^s_{i=1}\gamma_{k-1,i-1}},
\end{array}\end{equation}
empty sum being understood to be zero.

\begin{thm}
Suppose that $\varphi$ satisfies  \eqref{varphi-behavior} with $\theta=1$. Then for any  $z\in \mathbb{C}$ with $|\arg z|<\pi$ and $n\in \mathbb{N}^+$, we have
\begin{equation}\label{Stieltjes-expan-theta=1}
\begin{array}{l}
\displaystyle{  S_\varphi(z)= \sum^{n-1}_{s=0} \sum^s_{i=0} \sum^i_{j=0}\sum^{j+1}_{l=0}
\frac{(s+1)!  \alpha_{s, i}    C^j_{1, s, i}
\begin{pmatrix}
  j+1 \\
  l \\
\end{pmatrix}
\gamma_{j+1-l,s}}{j+1}
 z^{-s-1}\ln^l z
}
\\[.5cm]
\displaystyle{
~~~~~~~~~~     -
\sum^n_{s=1} (s-1)! c^*_s z^{-s} +\frac {(-1)^n}{z^n} \int^\infty_0\frac {\tau^n \varphi_n(\tau)}{\tau+z} d\tau,}
\end{array}
\end{equation}
where $c^*_s=\varphi_{s,s}(1)-\int_0^1 \varphi_{s-1,s-1}(t) dt$,
 $\varphi_{s,s}(t)=\frac{(-1)^s}{(s-1)!}\int^\infty_t(\tau-t)^{s-1} \varphi_s(\tau) d\tau$ for  $s\in \mathbb{N}^+$,
 and $\varphi_{0,0}(t)=\varphi(t)$.
\end{thm}\vskip .5cm

In view of the well-known formula of   Plemelj
\begin{equation}\label{Plemelj}
H^+_\varphi(x)=\frac 1 2 \lim_{\varepsilon\to 0^+} \int^\infty_0\left [
\frac {1}{t-x+i\varepsilon}- \frac {1}{t-x-i\varepsilon}\right ] \varphi(t) dt,
\end{equation}the results of the last two theorems also give the asymptotic expansions of the one-side Hilbert transform of $\varphi$.

\begin{thm}
Suppose that $\varphi$ satisfies  \eqref{varphi-behavior} with $\theta\in (0, 1)$. Then for any  $n\in \mathbb{N}^+$, we have
\begin{equation}\label{Hilbert-expan}
\begin{array}{l}
\displaystyle{ H^+_\varphi(x)= \sum^{n-1}_{s=0} \sum^s_{i=0} \sum^i_{j=0}
(-1)^{s+j} \pi \alpha_{s, i} C^j_{\theta, s, i}
\frac {d^j}{d\theta^j}\left [  (\theta)_s  \cot(\theta\pi) x^{-\theta} \right ] x^{-s}
}
\\[.5cm]
\displaystyle{
~~~~~~~~~~~     -
\sum^n_{s=1} (-1)^s (s-1)! c_s x^{-s} +\frac {1}{x^n}
\int^\infty_0{\!\!\!\!\!\!\!\!\!\!\!\! {}-{} \,\,\,\,\,  }\frac {\tau^n\varphi_n(\tau)}{\tau-x} d\tau,
}
\end{array}
\end{equation}
where $c_s=\frac{(-1)^s}{(s-1)!} M[\varphi; s]$.
\end{thm}\vskip .5cm

\begin{thm}
Suppose that $\varphi$ satisfies  \eqref{varphi-behavior} with $\theta=1$. Then for any   $n\in \mathbb{N}^+$, we have
\begin{equation}\label{Hilbert-expan-theta=1}
\begin{array}{l}
\displaystyle{  H^+_\varphi(x)= \sum^{n-1}_{s=0} \sum^s_{i=0} \sum^i_{j=0}\sum^{j+1}_{l=0}
\frac{(-1)^{s+1} (s+1)!  \alpha_{s, i}    C^j_{1, s, i}
\begin{pmatrix}
  j+1 \\
  l \\
\end{pmatrix}
\gamma_{j+1-l,s}}{2(j+1)}
}\\[.5cm]
\displaystyle{
~~~~~~~~~~~~~~~~~~~~~~~~~~~  \times \left [ (\ln x+i\pi)^l+(\ln x-i\pi)^l\right ]
x^{-s-1} }
\\[.5cm]
\displaystyle{
~~~~~~~~~~~     -
\sum^n_{s=1}(-1)^s (s-1)! c^*_s x^{-s} +\frac {1}{x^n}
\int^\infty_0{\!\!\!\!\!\!\!\!\!\!\!\! {}-{} \,\,\,\,\,  }
\frac {\tau^n \varphi_n(\tau)}{\tau-x} d\tau,}
\end{array}
\end{equation}
where $c^*_s=\varphi_{s,s}(1)-\int_0^1 \varphi_{s-1,s-1}(t) dt$,
 $\varphi_{s,s}(t)=\frac{(-1)^s}{(s-1)!}\int^\infty_t(\tau-t)^{s-1} \varphi_s(\tau) d\tau$ for  $s\in \mathbb{N}^+$, and $\varphi_{0,0}(t)=\varphi(t)$.
\end{thm}\vskip .5cm

To obtain the asymptotic behavior of the solution to equation \eqref{int-eq-WH}, one would have to make use of the results in the above four theorems.  For examples, see \cite{Li-Wong-2021}.

%

\end{document}